
\input amssym.def
\magnification=\magstep1

\hfuzz=11.1pt
\jot=14truept
\input psfig
%
\centerline{{\bf Axisymmetric Solutions of the Euler Equations }}
\centerline{{\bf for Sub-Square Polytropic Gases}}
	\vskip.75truein
\begingroup\baselineskip=12truept
\centerline{Yuxi Zheng}
\centerline{Department of Mathematics} 
\centerline{Indiana University}
\centerline{Bloomington, IN  47405}
		\bigskip
\centerline{Tong Zhang}
\centerline{Institute of Mathematics} 
\centerline{Academia Sinica}
\centerline{Beijing, China}
\endgroup\par

         \vskip 2truein

\footnote{}{ {\it AMS(MOS) Subject Classifications}. 
Primary: 35L65, 35L67, 65M06, 76N10;
Secondary: 65M99.}

\footnote{}{{\it Keywords}: Compressible, explicit solutions, 
hurricane, tornado, ordinary differential equations, 
Riemann problems, self-similar, swirl, two-dimensional, vacuum.}

\vfill\eject

\noindent 
{\it Abstract}.\quad 
We establish rigorously the existence of a three-parameter family of 
self-similar, globally bounded, and continuous weak solutions in two space 
dimensions to the compressible Euler equations with axisymmetry for 
$\gamma$-law polytropic gases with $1 \le \gamma < 2$. The initial data of 
these solutions have constant densities and outward-swirling velocities.
We use the axisymmetry and self-similarity assumptions to reduce the equations
to a system of three ordinary differential equations, from which we obtain 
detailed structures of solutions besides their existence. 
These solutions exhibit familiar structures seen in hurricanes and tornadoes. 
They all have finite local energy and vorticity with well-defined initial and 
boundary values.

\vskip.4truein
\noindent {\bf 1. Introduction.}\quad

  We are interested in finding some solutions to the initial value 
problem for the two-dimensional compressible Euler equations.
Our approach is to try to generalize to two dimensions some of the results
on Riemann problems for the one-dimensional compressible Euler equations.
One natural generalization is to consider initial data which 
consist of four constant states, or any finite number of 
constant states [13--15]. The well-known configurations of regular and Mach
reflections are special cases of such a generalization. However, no rigorous 
proof of existence of any nontrivial solutions to these problems has been 
established.

A more complete generalization is to consider initial data which 
depend only on the polar angle in the two-dimensional space of positions.
This generalization looks certainly more formidable, if one hopes to solve 
it all. However, it now contains a special three-parameter family of data, 
namely, the axisymmetric initial data, which allows us to
assume axisymmetry of the solutions and  reduce the initial 
value problem of the partial differential equations to a boundary value 
problem of a system of ordinary differential equations.  With this 
symmetry we were able to construct in [16] a two-parameter 
family of selfsimilar solutions corresponding to pure rotational initial data.
In paper [17], we gave rigorous proofs for these solutions and
constructed an additonal one-parameter family of solutions for the 
square polytropic cases. This additional parameter allows the initial
flows to swirl outward. In the current paper we take on the sub-square cases:
we construct rigorously a three-parameter family of solutions. 

Our main task here is to study the resulting boundary value problem 
of the non-autonomous system of three ordinary differential equations 
in the cases $1 \le \gamma < 2$. Singularity (stationary) points of the system 
consist of two-dimensional manifolds in the four-dimensional phase space. 
These singularity points 
correspond to surfaces of characteristics, where solutions may not be 
differentiable in the physical space and time. Linearizations at some 
stationary points yield no structure for the solutions due to the high 
order degeneracy of the system at these points. The existence of connecting 
orbits follows from high-order local analysis at stationary points 
and global analysis of various invariant regions of the system.
A typical global solution may consist of as many as three nontrivial 
connecting orbits chained together continuously. The complete construction 
of the three-parameter family of solutions is described in the conclusions 
section, section 8,  at the end of the paper.

We find that our solutions capture qualitatively some typical properties of 
swirling flows such as hurricanes and tornadoes. As an illustration, we plot 
a typical solution in Fig. 1.1--2.
These figures are placed out of the main flow of the text 
due to their large sizes.
See Fig. 8.1--2 for more precise description on the eye and wall regions. 
We refer the interested reader to paper [16] for explicit presentations of 
spiralling particle trajectories.

We point out that our initial data are restricted to the set of data with 
nonnegative radial velocities (swirling outward), in addition to axisymmetry
and radial symmetry.  Shock waves may arise if the initial radial velocities 
are allowed to be negative (swirling inward). Due to complications of the 
global dynamics of the ordinary differential equations,
we have only studied the special case $\gamma =2$ in a companion paper [17].
In this paper, we study the case $1 \le \gamma <2$. The case $\gamma > 2$ will
be considered in a forthcoming paper. For some partial results on the 
case $\gamma > 2$, we refer the reader to our paper [16].
Preliminary steps (Sect. 1--5) are reproduced here
from  [17] for readers'  convenience.

There is a great deal of related work, from which we mention only the 
most closely related. For general existence of weak solutions with axisymmetry 
for the 2-D compressible Euler equations outside a core region, we refer 
the reader to the recent work of Chen and Glimm [5].
For some explicit solutions of the compressible Euler equations
with spherical symmetry but without swirls, see Courant and Friedrichs [4].
For viscous swirling motions, we refer the reader to Bellamy-Knights [1],
Colonius, Lele, and Moin [3], Mack [9], Powell [11], and Serrin [12].

\vskip.4truein
\noindent {\bf 2. The problem.}\quad

We consider the two-dimensional compressible and  polytropic Euler equations
$$\eqalignno{
& \rho_t +(\rho u)_x + (\rho v)_y = 0,\cr
& (\rho u)_t + (\rho u^2 + p)_x + (\rho uv)_y = 0, &(2.1)\cr
& (\rho v)_t + (\rho uv)_x + (\rho v^2 + p)_y = 0, \cr
}$$
where $p =p(\rho)$ is a given increasing function of $\rho$. 
Global existence of
weak solutions to its initial value problem  is open. Attempts have
been made through considering special situations such as the diffraction
of a planar shock at a wedge or generalized Riemann problems with 
four different initial constant states. 

Here we consider a situation which involves swirling motions.
We impose axisymmetry to the system. That is, we assume that our solutions
$(u, v, \rho)$ have the property
$$\eqalignno{
\rho (t{,}r{,}\theta) &= \rho(t{,}r{,}0)\cr
\left( \matrix{
u(t{,}r{,}\theta)\cr
v(t{,}r{,}\theta)\cr
}\right) &= \left( \matrix{
\cos\theta & -\sin\theta\cr
\sin\theta & \cos\theta\cr
}\right) \left( \matrix{
u(t{,}r{,}0)\cr
v(t{,}r{,}0)\cr
}\right) &(2.2)\cr
}
$$
for all $t\ge 0$, $\theta \in \Bbb R$ and $r>0$, where $(r, \theta)$ is the
polar coordinates of the $(x, y)$ plane. With this symmetry, system (2.1)
can be reduced for continuous solutions to 
$$
\eqalignno{
& \rho_t +(\rho u)_r + {\rho u\over r} = 0,\cr
&  u_t + u u_r + {p_r\over \rho} - {v^2\over r} = 0, &(2.3)\cr
& v_t + uv_r + {uv\over r} = 0, \cr
}$$
where $\rho = \rho(t, r, 0)$, etc.  Note now that $u$ and $v$ in
(2.3) represent the radial and pure rotational (a.k.a. tangential )
velocities in the flow, respectively.

We limit ourselves to a Riemann-type of initial data; that is, 
we require the initial data to be independent of the radial variable
$r>0$:
$$
(\rho(0, r, \theta), u(0, r, \theta), v(0, r, \theta)) = (\rho_0(\theta),
u_0(\theta), v_0(\theta)). \eqno(2.4)
$$
We remark in passing that the initial value problem of system (2.1)
with the type of data in (2.4) seems to be more in the spirit to be called
the two-dimensional Riemann problem for (2.1).
It degenerates to the classical Riemann problem for the one-dimensional case,
 it is simple, and yet general enough to contain important waves such as
swirling motions as well as shock and rarefaction waves and 
slip lines (surfaces).

When the axisymmetry condition (2.2) is imposed onto (2.4), we find that
our data are limited to be
$$
\eqalignno{
u(0{,}r{,}\theta)    &=u_0\cos\theta - v_0 \sin\theta, \cr
v(0{,}r{,}\theta)    &=u_0\sin\theta + v_0\cos\theta, &(2.5)\cr
\rho(0{,}r{,}\theta) &= \rho_0, \cr
} 
$$
where $\rho_0 > 0, (u_0, v_0) \in \Bbb R^2$ are arbitrary constants.
Hence our data for system (2.3) are
$$
\rho(0, r, 0) = \rho_0,\quad u(0, r, 0) = u_0,\quad v(0, r, 0) = v_0. 
\eqno(2.6)
$$

Since the problem (2.3)(2.6) is invariant under self-similar transformations,
we look for self-similar solutions $(\rho, u, v)$ which depend only on 
$\xi = r/t$. We thus have the following boundary value problem of 
a system of ordinary differential equations:
$$
\rho_r = {\rho\over r} {\Theta\over \Delta},\quad u_r = {1\over r}
{\Sigma\over \Delta},\quad v_r = {uv\over r(r-u)}, \eqno(2.7)
$$
$$
\lim_{r\to+\infty} (\rho, u, v) = (\rho_0, u_0, v_0)\eqno(2.8)
$$
where 
$$
\eqalign{
\Delta &\equiv c^2 - (u-r)^2, \cr
\Theta &\equiv v^2-u(r-u),\cr
\Sigma &\equiv (r-u)\Theta - u\Delta = v^2(r-u) - uc^2, \cr
c      &\equiv \sqrt{p'(\rho)}, \cr
}
$$
and $r$ is used in place of $\xi$ for notational convenience.
We note that $c$ is the sound speed.

We will construct global continuous solutions or establish their existence for
problem (2.7--8) for any $\rho_0 >0$, $u_0 \ge 0$, $v_0 \in \Bbb R$ and
$p(\rho) = A_2\rho^\gamma$ where $1 \le \gamma < 2$ and $A_2 > 0$ is a 
constant. We shall be mainly concerned with the case $1 < \gamma < 2$ in
Sections 4--7, but we will mention the differences for $\gamma =1$ in
Section 8. We will also  assume $v_0 \ge 0$ because the case $v_0 <0$ can 
be transformed by 
$v \to -v$ to the case $v_0 > 0$. It can be verified that these
solutions are also global continuous solutions to the original 
Euler equations.

\vskip.4truein
\noindent {\bf 3. Far-field solutions.}\quad

We show that problem (2.7-8) has a local solution near $r = +\infty$ 
for  any datum $(\rho_0, u_0, v_0)$ with $\rho_0 > 0$. Let $s = {1\over r}$.
Then (2.7-8) can be written as 
$$
\eqalignno{
{d\rho\over ds}  & = {\rho[u(1-us)-sv^2]\over s^2c^2 - (1-us)^2} \cr
{du\over ds}  & = {suc^2 - v^2(1-us)\over s^2c^2 - (1-us)^2} 
                                           &(3.1) \cr
{dv\over ds}  & = -{uv\over 1-us}, \cr
}
$$
$$
(\rho, u, v)|_{s=0} = (\rho_0, u_0, v_0). \eqno(3.2)
$$
Problem (3.1-2) is a classically well-posed problem which has a 
unique local solution for any initial datum with $\rho_0 > 0$.

We find that $v=0$ and 
$$
sv^2-u(1-us) =0 \eqno(3.3)
$$
are invariant surfaces in the four-dimensional $(\rho, u, v, s)$ space.
It can be verified that 
$$
{d\over ds}[sv^2-u(1-us)] = {(1-3su)sc^2 + u(1-su)^2\over 
(1-su)[s^2c^2 - (1-su)^2]}[sv^2-u(1-us)].
$$

\bigskip
{\it Explicit solutions.}\quad
In the special case $u_0 = 0$, we find in [16] from the invariant surface (3.3)
a set of  explicit solutions near  $r = +\infty$ :
$$
\rho = \rho_0,\quad u = {v_0^2\over r}, \quad v = {v_0\over r}\sqrt{r^2-v_0^2},
\qquad r \ge r^* \eqno (3.4)
$$
where
$$
r^* \equiv {1\over 2}\left(\sqrt{c^2_0 + 4v_0^2} + c_0\right),\eqno (3.5)
$$
and 
$$
c_0 \equiv \sqrt{p'(\rho_0)}.
$$
We comment that the functions in (3.4) are actually defined for all 
$r\ge v_0$, but we cannot use them up to $v_0$ with absolute certainty 
since (3.1) has a singularity at the point $r=r^*$ on the curve (3.4).
In fact, we find that the position $r$, the radial velocity $u$, and the
sound speed $c$ at $r^*$ along (3.4) have the relation
$$
r^* = u(r^*) + c_0\ ,
$$
which is to say that $r^*$ is the radial characteristic speed; 
that is, distance that a small disturbance
generated from the origin at time zero can travel radially in time 
$t=1$. Thus, $r^*$ in (3.5) is a characteristic position.

We point out that the Mach number along any solution in (3.4) is a 
constant:
$$
M \equiv \sqrt{u^2+v^2}/c = v_0/c_0.
$$
For those who are interested in the pseudo Mach number defined in the
self-similar coordinates $(\xi, \ \eta)$ as 
$$
M_s \equiv \sqrt{(u-\xi)^2 + (v-\eta)^2}/c,
$$
we find that 
$$
M_s = \sqrt{r^2 - v_0^2}/{c_0}
$$
along a solution of (3.4). At the end point $r=r^*$, we find
$$
M_s(r^*) = \left(\sqrt{M_0^2+1/4}+1/4\right)^{1/2} > 1,
$$
where 
$$M_0 \equiv v_0/{c_0}$$ 
is the Mach number for the initial data in this
case. Pseudo Mach number determines the type (hyperbolic or other) of the 
2-D Euler equations in the self-similar plane $(\xi, \ \eta)$: hyperbolic
if $M_s > 1$, other types if $M\le 1$, see [15] for example.
So the solutions (3.4) are pseudo supersonic and in hyperbolic regions.

\vskip.4truein
\noindent {\bf 4. Intermediate field equations.}

We can simplify system (3.1) by assuming the polytropic relation
$$
p(\rho) = A_2\rho^\gamma \eqno (4.1)
$$
for some $A_2 > 0$ and $\gamma > 1$, and introducing the variables
$$
I = su, \quad J = sv, \quad K = sc. \eqno(4.2)
$$
Then system (3.1) can be written into the form
$$
\eqalignno{
s{dI\over ds} & = {2IK^2 - (1-I)[J^2+I(1-I)]\over K^2 - (1-I)^2} \cr
s{dJ\over ds} & = J{1-2I\over 1-I} &(4.3)\cr
s{dK\over ds} & ={K\over 2}{2K^2 - 2(1-I)^2 - (\gamma-1)[J^2 - I(1-I)]\over
K^2 - (1-I)^2}. \cr
}
$$
Corresponding to the  initial data (3.2), we shall look for solutions
of (4.3) with the following initial condition:
$$
(I, J, K) \sim s(u_0, v_0, c_0) \eqno(4.4)
$$
as $s \to 0+$.
We note that (4.3) is now autonomous for $I, J, K$ with respect to
the new variable $s' = \ln s$. 

The invariant surfaces of (4.3) are the surface $J = 0,$ the surface
$K =0$, and the surface 
$$
H \equiv J^2 - I(1-I) = 0 \eqno(4.5)
$$
since 
$$
s{d\over ds}H = {(1-2I)[ 2K^2 - 3(1-I)^2]\over (1-I)[K^2 - (1-I)^2]}H.
$$
The explicit solutions (3.4) are all in the invariant surface $H=0$.
We postpone the description of these solutions in the variables 
$(I, J, K)$ to Sect. 6 when we have a more complete picture of all solutions.

\smallskip
{\it Scaling symmetry}.\quad 
System (4.3) is invariant under the coordinate transformation
$ s \to \alpha s$ for any constant $\alpha > 0$. In particular, we  can
take $\alpha = 1/{c_0}$. Thus we may assume that $\rho_0 > 0$
 is such that $c_0 = \sqrt{p'(\rho_0)} = 1$. 
Hence the structure of any solution
of problem (4.3-4) will depend only on the 
dimensionless quantities  $u_0/{c_0}$ and $v_0/{c_0}$.

\vskip.4truein

\noindent {\bf 5.  Solutions without swirls.}
		\medskip
Let us first determine the distribution of integral curves on the
invariant surface $J=0$.

Assume $v_0 = 0$.  We look for solutions to problem (4.3-4) with $J=0$.  Hence
we have a subsystem for $(I,K)$:
$$\cases{
s {dI\over ds} = I {2K^2 - (1-I)^2\over K^2 - (1-I)^2} \cr
\noalign{\vskip10truept}
s {dK\over ds} = K {K^2 - (1-I)^2 + {\gamma-1\over 2} I(1-I)\over K^2 -
(1-I)^2} \cr
}\eqno\eqalign{
&(5.1) \cr
&(5.2)\cr
}$$
Introducing a new parameter $\tau$, we can rewrite (5.1-2) as
$$\cases{
{dI\over d\tau} = I\left[ (1-I)^2 - 2K^2\right] \cr
		\noalign{\vskip12truept}
{dK\over d\tau} = K\left[ (1-I)^2 - K^2 - {\gamma-1\over 2} I(1-I)\right]\cr
		\noalign{\vskip12truept}
{ds\over d\tau} = s\left[ (1-I)^2 - K^2 \right]. \cr
}\eqno\eqalign{
&(5.3) \cr
&(5.4) \cr
&(5.5) \cr
}$$
Note that (5.3-4) form an autonomous subsystem.  If $u_0$ also vanishes, then
we have a trivial solution $\rho = \rho_0,\ u = v = 0$.

Suppose $u_0 > 0$.  It can be verified that our far-field solutions starting at
$s = 0+$ enter the region $\Omega\subset \Bbb{R}^2$ in the $(I,K)$ phase
space given by
$$\Omega:\cases{
I>0,\ K>0, \cr
		\noalign{\vskip6truept}
a\equiv (1-I)^2 - K^2 - {\gamma-1\over 2}I(1-I)>0 & in \quad$0<I\leq 
{1\over \gamma}$\cr
		\noalign{\vskip6truept}
b\equiv (1-I)^2 - 2K^2 > 0 & in \quad${1\over \gamma} \leq I < 1$\cr
}$$
See Figure 5.1.

\medskip
\midinsert\centerline{\psfig{figure=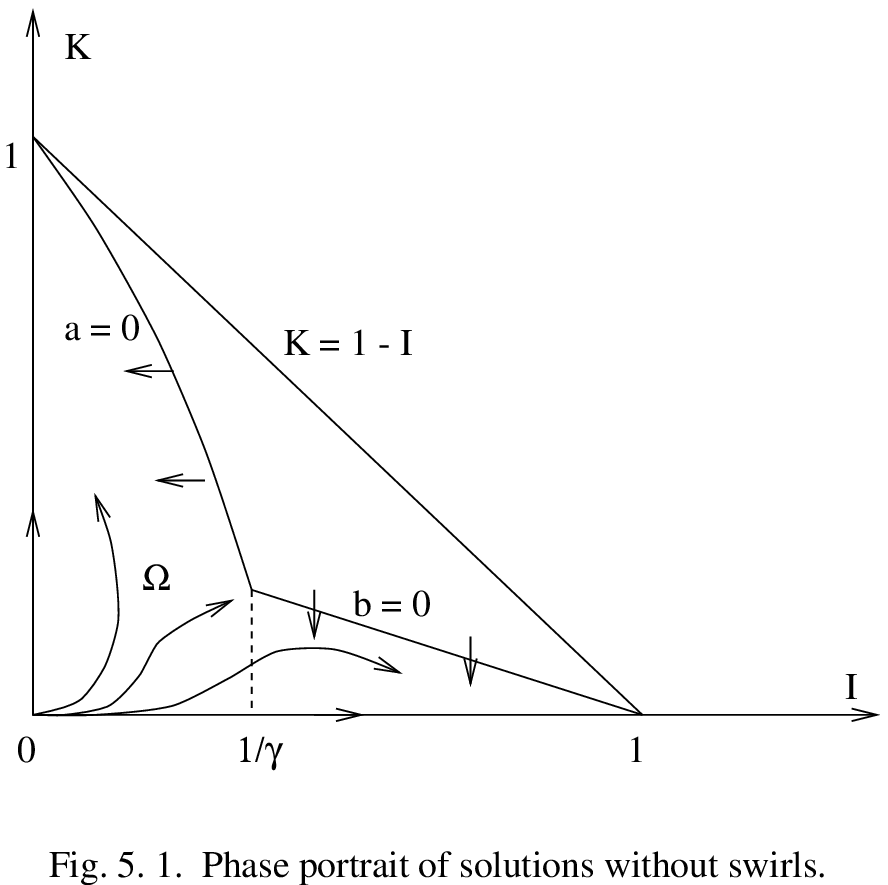,height=3.5truein}}
\endinsert
\medskip\medskip \medskip\medskip\medskip

We show that solutions starting in the closure
$\overline{\Omega}$ will not leave $\overline{\Omega}$ as $s$ increases. 
Notice first that $s>0$ is an increasing function of $\tau$ in $\Omega$ by
equation (5.5) so we will show that solutions of the two equations (5.3-4)
starting in the closure $\overline{\Omega}$ will not leave $\overline{\Omega}$
as $\tau$ increases.  The stationary points of (5.3-4) in $\overline{\Omega}$
are the points $(I,K) = (0,0), (0,1), (1,0)$, and 
$\left( {1\over \gamma}, {1\over
\sqrt{2}}\left( 1 - {1\over \gamma}\right)\right)$.  The axis $K=0$ in $0<I<1$,
and the axis $I=0$ in $0<K<1$ are trivial solutions.  On the boundary 
$b = 0$, i.e., 
$$
K = {1\over \sqrt{2}}(1-I),\quad {1\over \gamma}< I < 1, \eqno(5.6)
$$
we find that ${dI\over d\tau} = 0,\ {dK\over d\tau} < 0$.  So solutions enter
$\Omega$ on (5.6).  Finally on the boundary $a=0$, i.e., 
$$
K^2 = (1-I)^2 - {\gamma-1\over 2}I(1-I),\qquad  0<I<{1\over \gamma}, \eqno(5.7)
$$
we find that ${dI\over d\tau}<0$, and ${dK\over d\tau} = 0$.  
So solutions enter $\Omega$ on (5.7) also.

We conclude without showing further details that some solutions in $\Omega$
will go to the point $(1,0)$, while others go to the point $(0,1)$, 
with exactly one integral curve (the heteroclinic orbit)  going to the point 
$\left( {1\over \gamma}, {1\over \sqrt{2}}
\left(1-{1\over\gamma}\right)\right)$.

\medskip
\midinsert\centerline{\psfig{figure=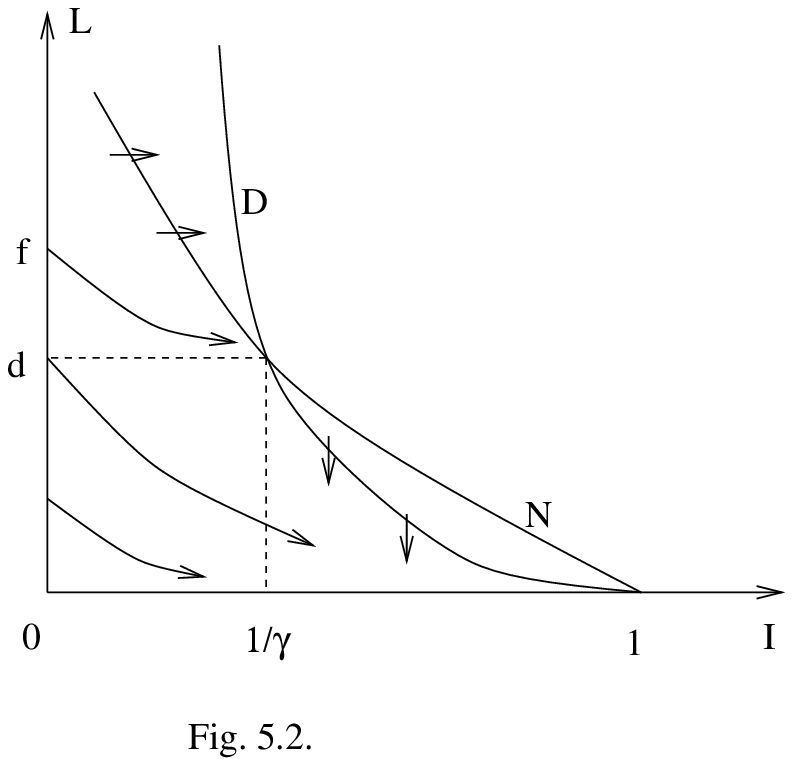,height=3truein}}
\endinsert\medskip\medskip\medskip\medskip\medskip

{\it 5.1. The transitional solutions.} 
To determine what data $(u_0,\rho_0)$ yield the transitional solution leading
to the stationary point $\left( {1\over \gamma},{1\over \sqrt{2}} \left(
1-{1\over
\gamma}\right)\right)$, we first eliminate $\tau$ from (5.3-4) by division and
introduce $L \equiv ( K/I)^2$ to find
$$
{dL\over dI} = -L {(\gamma-1)(1-I)-2IL\over (1-I)^2 - 2I^2L},\quad 0<I<1
\eqno(5.8)
$$
$$
L = M^{-2}_0 \ \hbox{at}\ I=0 \eqno(5.9)
$$
where 
$$
M_0 \equiv u_0/{c_0}
$$ 
denotes the Mach number of the initial
states $(u_0,0,\rho_0)$.  Problem (5.8-9) is well-posed for any $M_0>0$.  The
transitional solution goes from the point 
$(I,L) = (0,M^{-2}_0)$ to $\left( {1\over
\gamma},{1\over 2}(\gamma-1)^2\right)$.  However, there does not seem to 
have an explicit formula for the value $M_0$ which yield the transitional 
solution.  We give an estimate instead.  For convenience we let 
$M_h(\gamma)$ denote the value of the initial Mach number for the 
transitional (critical, heteroclinic) solution.  It can be seen that $M_0 =
{\sqrt{2}\over \gamma-1}$ is an upper bound for the transitional Mach number
$M_h(\gamma)$.  In fact, any solution $L(I)$ of (5.8-9) with datum $M_0\geq
{\sqrt{2}\over \gamma-1}$ starts as a decreasing function of $I\geq 0$ 
till $I =
{1\over \gamma}$.  In the interval $I\in \left[ {1\over \gamma},1\right]$, the
solution remains below the two curves on which the numerator and demoninator of
the right-hand side of (5.8) vanish respectively, therefore remains decreasing
until the final stationary point $(I,L) = (1,0)$, see Figure 5.2, where $d=
(\gamma -1)^2/2$, $f=M_h^{-2}(\gamma)$, and $N$ and $D$ are where the numerator
and denominator of the right-hand side of (5.8) vanish respectively.

This transitional solution yields a one-parameter family of smooth solutions in
terms of $(r,u,v,\rho)$.  We see from equation (5.5) that $\ln s$ approaches
infinity as the solutions approach the point $(I,K) = \left( {1\over \gamma},
{1\over \sqrt{2}}\left(1-{1\over \gamma}\right)\right)$, since $(1-I)^2 - K^2
\ne 0$ at the point.  So $s\to +\infty$ and thus $r\to 0+$.  Also the solutions
have the asymptotics
$$
u(r) = {1\over \gamma}r,\quad\ c(r) = {1\over \sqrt{2}}\left(1-{1\over
\gamma}\right)r
$$
as $r\to 0+$.  These global transitional solutions are similar and one is
sketched in Figure 5.3.

\medskip
\midinsert\centerline{\psfig{figure=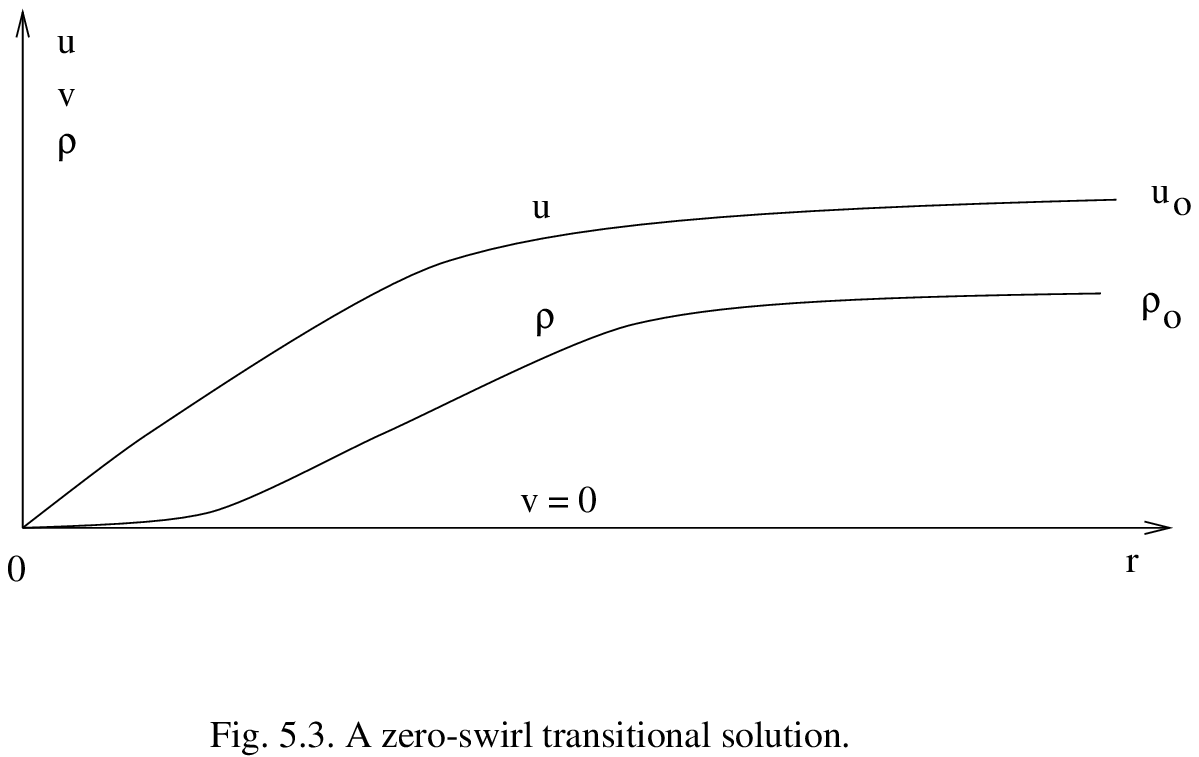,height=3.0truein}}
\endinsert
\medskip\medskip\medskip\medskip\medskip\medskip\medskip

The Mach number at $r=0$ of the transitional solution is 
$M = \sqrt{2}/(\gamma -1)$. The corresponding pseudo Mach number is 
$M_s = \sqrt{2}$.

{\it 5.2. Finiteness of the parameter $s$.}
We show next that the parameter $s$ approaches finite numbers when solutions of
(5.3-5) approach either the points $(I,K) = (1,0)$ or $(0,1)$.

We can linearize the two equations (5.3-4) at $(I,K) = (0,1)$ to find
$$\cases{
{dI\over d\tau} = -I \cr
	\noalign{\vskip10truept}
{d(K-1)\over d\tau} = -{\gamma+3\over 2} I-2(K-1). \cr
}\eqno(5.10)
$$
The eigenvlaues of (5.10) are $\lambda_1 = -1, \lambda_2 = -2$.  
So solutions of
(5.3-4) near $(0,1)$ approach $(0,1)$ exponentially as $\tau \to +\infty$. 
From equation (5.5), we find
$$
\ln {s\over s_0} = \int^\tau_{\tau_0} \left[ (1-I)^2 - K^2\right]d\tau 
$$
for some constants $s_0>0$ and $\tau_0$.  So $s$ approaches a finite number as
$\tau\to\infty$ since $(1-I)^2-K^2$ approaches zero exponentially.

Linearization at the point $(I,K) = (1,0)$ of the two equations (5.3-4) yield
the trivial system of zero right-hand side.  We need a different approach.  We
show that solutions near $(1,0)$ will enter $(1,0)$ in the sector bounded by
$K=0$ and the line
$$
K = \alpha(1-I)\eqno(5.11)
$$
for some $0<\alpha < {1\over \sqrt{2}}$, such that $\alpha^2 + {\gamma-3\over
4}>0$.  See Figure 5.4.

\medskip
\midinsert\centerline{\psfig{figure=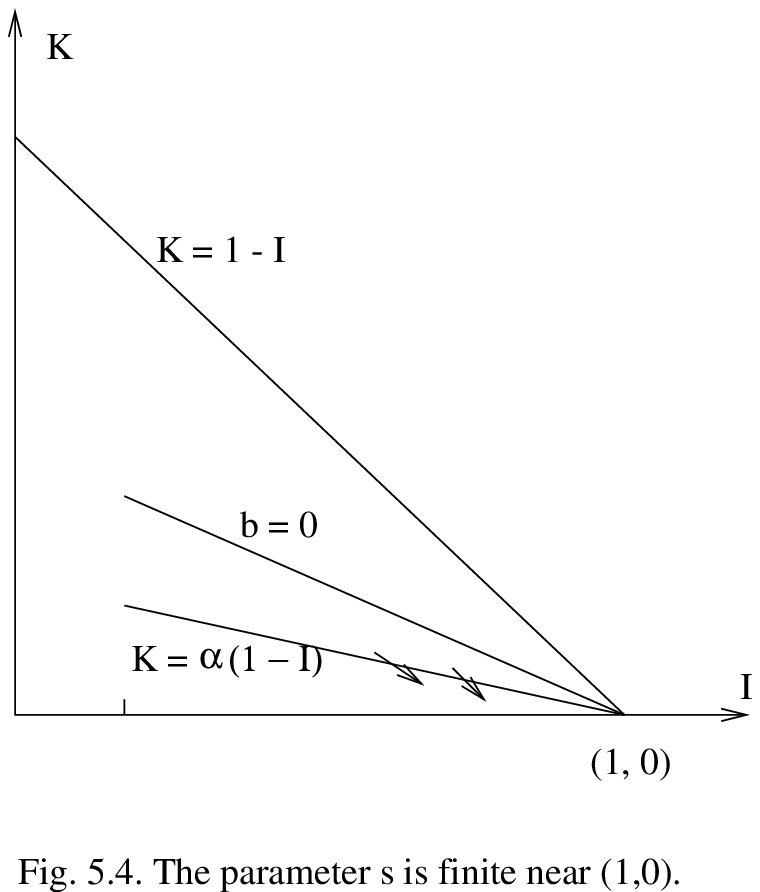,height=3.5truein}}
\endinsert
\medskip\medskip\medskip\medskip\medskip\medskip\medskip

In fact, a vector in the normal direction of (5.11) is $(\alpha,1)$.  We
calculate the inner product of the vector field of (5.3-4) with the direction
$(\alpha,1)$ to find
$$
{d\over d\tau}(I,K)\cdot (\alpha,1) = -\alpha \left[ 2I\left( \alpha^2 +
{\gamma-3\over 4}\right) - (1-\alpha^2)(1-I)\right](1-I)^2. \eqno(5.12)
$$
The expression (5.12) is negative if $I$ is close to 1 and $\alpha$ is such
that $\alpha^2 + {\gamma-3\over 4}>0$.  So we conclude that every solution that
ends at $(1,0)$ will be such that 
$$
K<\alpha(1-I) \eqno(5.13)
$$
near $I = 1$ for some $\alpha < {1\over \sqrt{2}}$, since $\gamma > 1$.  

Now we look at the equation (5.1) and use (5.13) to find
$$
s {dI\over ds} > I {1-2\alpha^2\over 1-\alpha^2}
$$
when $I$ is close to $1$.  Thus
$$
s{dI\over ds} > \hbox{a positive constant} 
$$
when $I$ is near $1$.  Therefore $s$ is finite at the end $(1,0)$ since $I$ is
bounded and the geometric integral $\int^{+\infty}_1 {1\over s} ds$ diverges.

{\it 5.3. The global construction.}
We are now ready to construct global solutions for (5.1-2).  For each solution
ending at $(I,K) = (0,1)$, we continue the solution by the constant state
$(u,v,\rho) = (0,0,\rho^{**})$ where $\rho^{**}$ is the value of
$\rho$ at the ending point.  These are continuous extensions.  The relation
between the density $\rho^{**}$ and the ending value $s^{**}$ is
$$p'(\rho^{**}) = s^{**-2}.$$
Let 
$$ c^{**} = \sqrt{p'(\rho^{**})},$$
then 
$$c^{**} =r^{**}$$
when $r^{**}$ is the radius of the circle of the constant state in the physical
plane $t=1$. So the core region has a constant density which expands at its
sound speed. The Mach and pseudo Mach numbers at the edge of the core region
are $0$ and $1$ respectively. In the core, the pseudo Mach number 
$M_s = r/c^{**} < 1$. 
 We therefore have constructed global solutions in this case.

For each solution ending at the point $(I,K) = (1,0)$, we continue the solution
by the vacuum state $\rho = 0$.  We do not need to specify the functions $u$ or
$v$ in the vacuum since the Euler equations have $\rho$ as a factor in every
term.  Each vacuum occupies a circular region of radius $r^{**}$ determined by
$r^{**} = {1\over s^{**}} = u^{**}$ in the physical plane at $t=1$, where
$u^{**}$ is the terminal radial velocity of the fluid at the edge of the
vacuum.  Asymptotically, we find (see also [16])
$$
\eqalign{
& u-r = -{2(\gamma -1)\over \gamma +1}(r - r^{**})(1+o(1)), \cr
& c = {(\gamma -1)\sqrt{3-\gamma}\over \gamma + 1}(r - r^{**})(1+o(1)) \cr
}$$
at the edge for $\gamma <3$.  For $\gamma =3$, the asymptotics is given by
$$
\eqalign{
& u-r = -(r - r^{**})(1+o(1)), \cr
& c = {1\over 2}(r - r^{**})\ln^{-1/2}{1\over r - r^{**}} (1+o(1)). \cr
}$$
And for $\gamma > 3$, it is given by
$$
\eqalign{
& u-r = -(r - r^{**})(1+o(1)), \cr
& c = K_0(r - r^{**})^{{\gamma-1\over 2}}(1+o(1)),\quad K_0 >0. \cr
}$$
From these we find that the Mach number $M=\infty$ for all $\gamma > 1$ and
the pseudo Mach number $M_s = {2\over \sqrt{3-\gamma}}$  for $\gamma < 3$,
and $M_s = \infty$ for $\gamma \geq 3$ at the edge.

In all, we have constructed global solutions to the reduced system
(5.1-2), the special case of (4.3-4) with zero swirl.
		\bigskip
\noindent {\bf 6.  General solutions in the intermediate field.}
		\medskip
Now consider the case $v_0>0$ and $u_0\geq 0$ for system (4.3).  Let
$\Omega_3\subset \Bbb{R}^3$ be the set of points $(I,J,K)$ 
satisfying $0<I<1,\ J>0,\ K>0$,
$$
H = J^2 - I(1-I)<0, 
$$
$$
B  \equiv (1-I)[J^2 + I(1-I)] - 2IK^2>0\ \  \hbox{if}\ \ 
 {1\over \gamma}\leq I < 1,
$$
$$
A \equiv 2(1-I)^2 + (\gamma-1)[J^2 - I(1-I)] - 2K^2>0\ \  \hbox{if}\ \ 
 0<I<{1\over
\gamma}.
$$
See Figure 6.1.  Note that the surface $A=0$ intersects the
coordinate plane $K=0$ on the line
$$
J^2 = {\gamma + 1\over \gamma -1}(1-I)(I - {2\over \gamma+1})
$$
which lies inside $H<0$. 
It can be verified that all far-field solutions with $u_0\ne
0,\ v_0>0, \rho_0>0$ enter the region $\Omega_3$ in $s>0$.  Far-field solutions
with $u_0=0$ enter the side $H = 0$ of $\Omega_3$, see Sect. 4.  
We omit these tedious verifications. 

Also let $\Omega_{33}\subset \Bbb{R}^3$ be the set of points $(I,J,K)$ 
satisfying $0<I<1,\ J>0,\ 0<K<1-I$, and $H<0$. See Fig.6.1.

We find that it is convenient to introduce a new variable $\tau$, as in section
5, to write the system (4.3) in the following form
$$\eqalignno{
{dI\over d\tau} & = (1-I)B & (6.1) \cr
{dJ\over d\tau} & = J(1-2I)[(1-I)^2 - K^2] & (6.2) \cr
{dK\over d\tau} & = {1\over 2}K(1-I)A & (6.3)\cr
{ds\over d\tau} & = s(1-I)[(1-I)^2 - K^2] & (6.4) \cr
}$$
This is an autonomous system for $(I,J,K,s)$, and the first three equations
(6.1-3) form an autonomous sub-system for $(I,J,K)$.

The stationary points of the system (6.1-3) contained in the closure
$\overline{\Omega}_{33}$ are found to be given by the edge
$$
K = 1-I,\ J^2 = I(1-I),\ \forall\, I\in [0,1] \eqno(6.5)
$$
and
$$
(I,J,K) = \left( {1\over \gamma},\ 0, \ {1\over \sqrt{2}}\left( 1-{1\over
\gamma}\right)\right) \eqno(6.6)
$$
in the case $\gamma \ne 2$.  For the case $\gamma = 2$, the stationary points
of (6.1-3) are given by edge (6.5) and the curve
$$
\cases{
I = {1\over 2}, \cr
		\noalign{\vskip6truept}
K^2 = {1\over 2}\left( J^2+{1\over 4}\right), \cr
}\quad 0\leq J < {1\over 2}.\eqno(6.7)
$$
which is also the intersection of $A=0$ with $B=0$ in this case.  
Hence there is no stationary point in the open region $\Omega_{33}$ when 
$1<\gamma<2$; all the stationary points are on the
boundary of $\Omega_{33}$ in this case.

\medskip
\midinsert\centerline{\psfig{figure=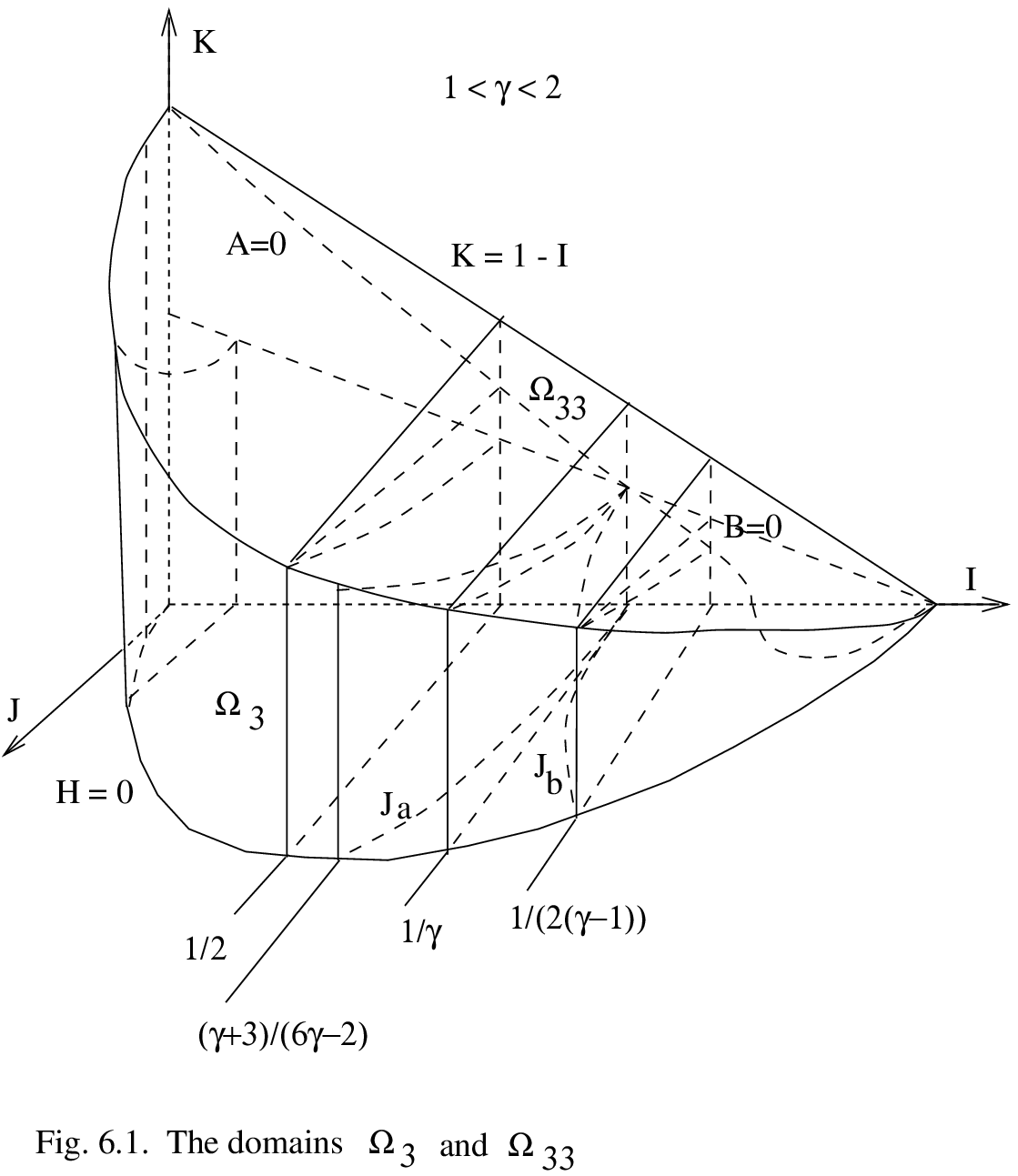,height=5.125truein}}
\endinsert
\medskip\medskip\medskip\medskip\medskip\medskip\medskip\medskip\medskip

We claim that solutions inside $\Omega_{33}$ do not leave $\Omega_{33}$ 
from its sides (excluding possibly edges, corners, or stationary points) as s
increases for all $\gamma > 1$.  
Note that the sides of $\Omega_{33}$ in the surfaces $K = 0$, or $J = 0$ or
$H = 0$ are invariant regions.  We need only to prove that no solution leaves
 $\Omega_{33}$ from the top side $K = 1- I$.
The top side has an outward normal
$$
\buildrel\rightharpoonup\over {n} = (1, 0, 1)
$$
arranged in the coordinate order $(I,J,K)$. 
We calculate the inner product of the normal
$\buildrel\rightharpoonup\over{n}$ with the tangent vector of an integral
curve of (6.1-3) on the top surface to yield
$$
\buildrel\rightharpoonup\over{n} \cdot {d\over d\tau}(I,J,K) =
(1-I)(B+{1\over 2}KA) < 0 
$$
which proves the claim, since ${ds\over d\tau}>0$ in $\Omega_{33}$. 

We analyze the local structure of solutions at the stationary point (6.6).
Let
$$
{\overline I} = I -\gamma, \quad {\overline J} = J,\quad {\overline K} = K -
{\gamma -1\over \sqrt{2}\gamma}.
$$
To first order system (6.1--3) can be reduced to 
$$\eqalignno{
{d{\overline I}\over d\tau} & = -{2(\gamma -1)^2\over \gamma^3}(
{\overline I}+\sqrt{2}\ {\overline K}) &  \cr
{d{\overline J}\over d\tau} & = -{(2-\gamma)(\gamma - 1)^2\over 2\gamma^3}
{\overline J} & \cr
{d{\overline K}\over d\tau} & = -{(\gamma -1)^3\over \sqrt{2}\gamma^3}
({\gamma + 2\over 2}{\overline I} + \sqrt{2}\ {\overline K}). & \cr
}$$
Its three eigenvalues are found to be
$$
\eqalignno{
\lambda_{\pm} & = -{(\gamma -1)^2\over 2\gamma ^3}\left[\gamma +1 \pm 
\sqrt{(\gamma+1)^2 + 4\gamma(\gamma-1)}\right] \cr
\lambda_2     & = -{(2-\gamma)(\gamma - 1)^2\over 2\gamma^3}. \cr
}$$
Corresponding eigenvectors are found to be
$$
\eqalignno{
\buildrel\rightharpoonup\over{v}_\pm & = \left(4\sqrt{2},\ 0,\ \gamma -3 \pm
\sqrt{(\gamma+1)^2 + 4\gamma(\gamma-1)}\right) \cr
\buildrel\rightharpoonup\over{v}_2 & =  (0,\ 1,\ 0). \cr
}$$
It follows that $\lambda_+ <0, \lambda_2 <0$, and $\lambda_- > 0$ 
for $1< \gamma < 2$. Hence the stationary point (6.6) is hyperbolic.
The local structure is depicted in Fig. 6.2.

To determine monotonicity of solutions, we need to 
deal with the surfaces $A=0$ and $B=0$ inside $\Omega_{33}$.
Consider first the surface
$$
B = I(1-I)^2 + (1-I)J^2 - 2IK^2 = 0\quad \hbox{in}\ \ 0 < I < 1.
\eqno(6.8)
$$
An outward normal is given by in the coordinate order $(I,J,K)$
$$
\left(2K^2 - (1-I)^2 + 2I(1-I) + J^2,\ -2(1-I)J,\ 4IK\right) \equiv
\buildrel\rightharpoonup\over {n}_B.
$$
We calculate the inner product of the normal
$\buildrel\rightharpoonup\over{n}_B$ with the tangent vector of an integral
curve of (6.1-3) on the surface (6.8) to yield
$$\buildrel\rightharpoonup\over{n}_B \cdot {d\over d\tau}(I,J,K) =
(1-I)^2HJ_b,
$$
where 
$$
J_b \equiv (\gamma I-1)(1-I) -(2- \gamma)J^2.
$$
See figure 6.1 for the position of $J_b = 0.$ 
So integral curves of (6.1--3) go
 from $B>0$ to $B<0$ on the surface $B=0$ with $J_b < 0$, but reverse their
directions on  the surface $B=0$ with $J_b > 0$ as $\tau$ increases.

\medskip
\midinsert\centerline{\psfig{figure=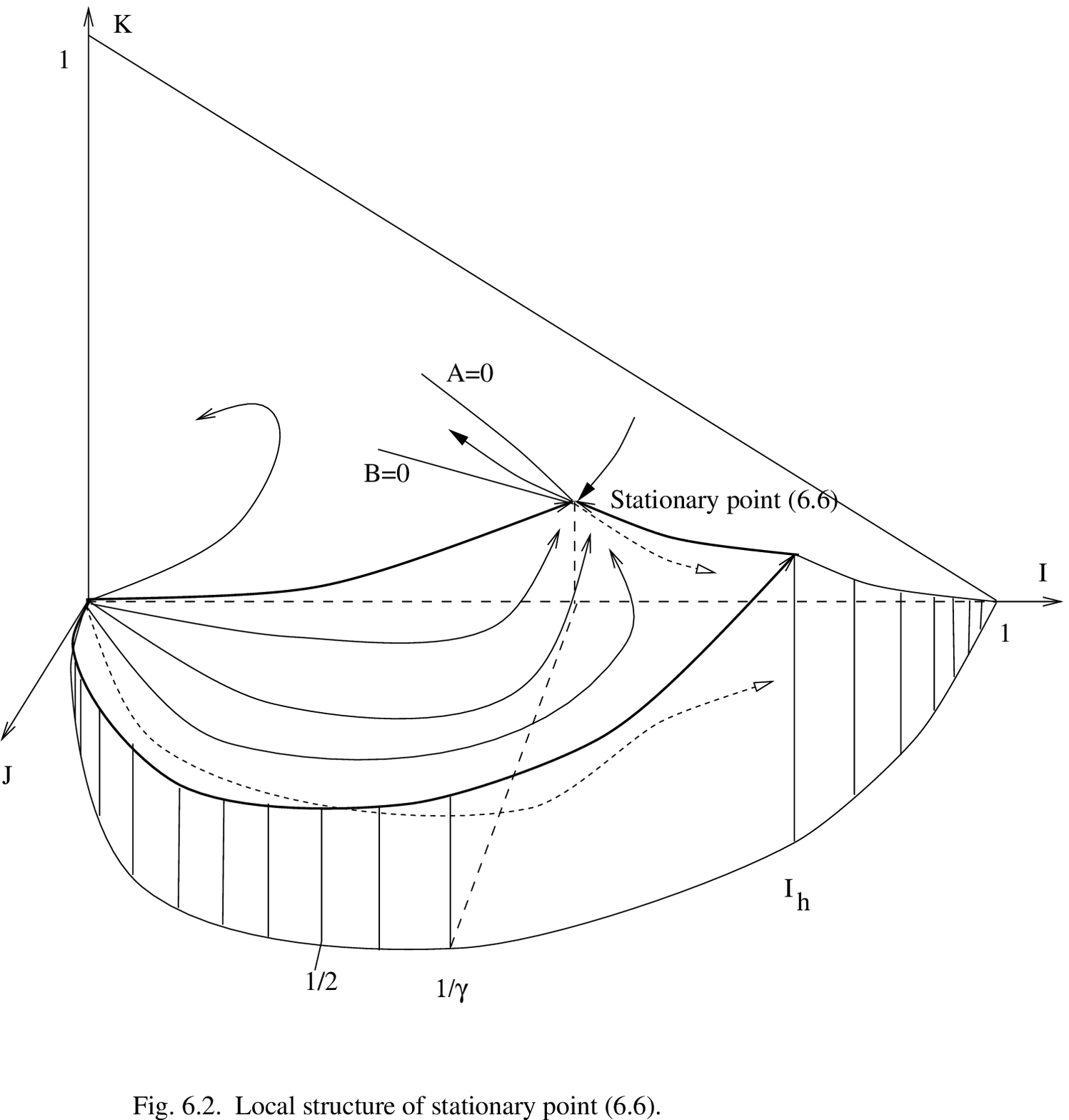,height=6.25truein}}
\endinsert
\medskip\medskip\medskip\medskip\medskip\medskip\medskip\medskip\medskip

\ 

\ 

Now consider the surface
$$
A = 2(1-I)^2 + (\gamma-1)[J^2-I(1-I)]-2K^2 = 0,\quad 0<I<1.
$$
An outward normal is given by
$$
\left(4(1-I)+(\gamma-1)(1-2I),\ -2(\gamma-1)J,\ 4K\right) \equiv
\buildrel\rightharpoonup\over{n}_A
$$
in the order $(I,J,K)$.  We similarly calculate the inner product of the normal
$\buildrel\rightharpoonup\over{n}_A$ with the tangent of the integral curves on
this surface to yield 
$$
\buildrel \rightharpoonup\over{n}_A\cdot {d\over d\tau} (I,J,K)  
= HJ_a,
$$
where 
$$
J_a \equiv [\gamma+3 - 2(\gamma+1)I] (1-I)(1-\gamma I) + (\gamma-1)^2J^2(1-2I).
$$
See figure 6.1 for the position of $J_a = 0.$  In particular the 
cylindrical surface $J_a = 0$ intersects the $(I,\ 0, \ K)$-plane 
at $I= 1/\gamma$ and the surfce $H=0$ 
at $I= {\gamma+3\over 2(3\gamma -1)}$.
Therefore integral curves of (6.1--3) go from $A>0$ to $A<0$ on
the surface $A=0$ with $J_a < 0$, but reverse their directions
on the surface with $J_a > 0$  as $\tau$ increases. The portion of 
the surface $A=0$ below $K=0$ will not be used.

\medskip
\midinsert\centerline{\psfig{figure=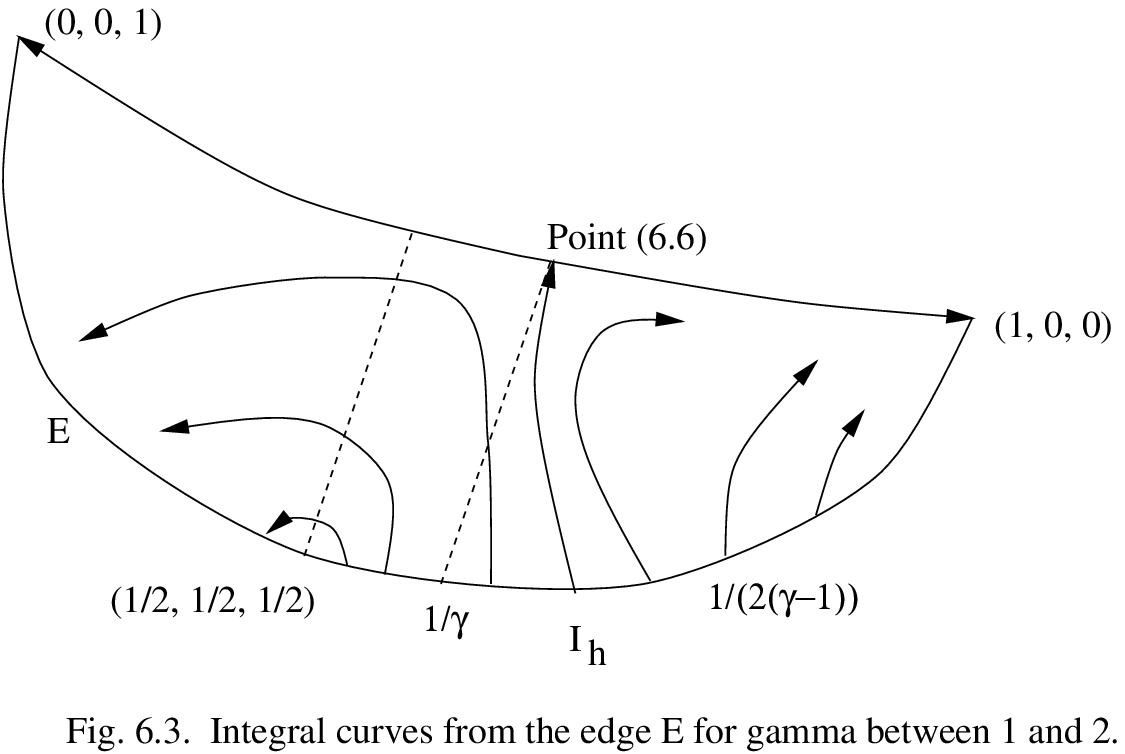,height=3.0truein}}
\endinsert
\medskip\medskip\medskip\medskip\medskip\medskip\medskip

We need to study the local structure of solutions near stationary points given
in formula (6.5).  It is helpful to first describe solutions of (6.1-4) with
data $u_0=0,\ v_0>0,\ \rho_0>0$.  Our far-field solutions (3.4-5) can be
written in terms of $I,J,K$ as follows:
$$\cases{
I  = s^2v^2_0 \cr
		\noalign{\vskip6truept}
J = sv_0\sqrt{1-s^2v^2_0}\cr
		\noalign{\vskip6truept}
K = sc_0\cr
}\eqno(6.9)
$$
valid for $s \in [0,s^*]$ where $s^* = 1/r^*$.  These solutions are all in the
surface $H = 0$ and they all start from the origin $(I,J,K) = (0,0,0)$ and end
at points of the stationary edge given by (6.5).  As the initial Mach number
$M_0 = |v_0|/{c_0}$ varies in $(0,\infty)$, the ending points
of the solutions cover all the interior points of the stationary edge (6.5)
exactly once.  

Now we study the structure of solutions near the stationary edge (6.5).  More
precisely, let us use $E$ to denote the open set
$$
E \equiv \{ (I,J,K) = (\alpha,\ \sqrt{\alpha(1-\alpha)},\ 1-\alpha) 
\quad | \quad \forall \alpha \in (0,1)\} \eqno(6.10)
$$
which contains the interior points of the curve (6.5).  
At any point of $E$, the linear part of the right-hand side of (6.1--3)
is found to be given by the matrix
$$
-(1-\alpha)\left[\matrix{
(1-\alpha)(1+2\alpha) 
& -2(1-\alpha)J_\alpha
& 4\alpha(1-\alpha) \cr
-2(2\alpha -1)J_\alpha
& 0 
& -2(2\alpha -1)J_\alpha \cr
(1-\alpha)\left[{\gamma-1\over 2}(1-2\alpha)+2(1-\alpha)\right] 
& -(1-\alpha)(\gamma -1)J_\alpha
& 2(1-\alpha)^2 \cr
}\right] 
$$
where $J_\alpha \equiv \sqrt{\alpha(1-\alpha)}$.
The eigenvalues and eigenvectors of this matrix can be found as follows.
Since the whole line $E$ is stationary, we expect that $\lambda_1 = 0$
 is an eigenvalue,
and the tangent vectors of $E$ are associated eigenvectors. The
solutions of (6.9) offer another set of eigenvectors
which can further simplify the calculation. We find the other two
eigenvalues to be 
$$
\lambda_2 = -2(1-\alpha)^2(1+\alpha)< 0,\quad 
\lambda_3 = -(1-\alpha)^2(1-2\alpha)
$$
with associated eigenvectors
$$
\eqalignno{
\buildrel \rightharpoonup\over{v}_2  &= (2\sqrt{\alpha},\
{1-2\alpha\over \sqrt{1-\alpha}},\  {1-\alpha\over \sqrt{\alpha}}) \cr
\buildrel \rightharpoonup\over{n}_3  & \equiv 
( 2(1-\alpha)[2(\gamma-1)\alpha-1],\ \   
-2(\gamma+1)\sqrt{\alpha(1-\alpha)}(2\alpha-1), &(6.11)\cr
&\qquad \qquad -(1-\alpha)[2(3\gamma-1)\alpha - (\gamma+3)]). \cr
}$$
Along this direction (6.11), we calculate the following information:
$$\eqalignno{
\buildrel\rightharpoonup\over{n}_3 \cdot
\buildrel\rightharpoonup\over{n}_A\big|_E & = 
2(1-\alpha)(2\alpha-1)[2(3\gamma-1)\alpha -  (\gamma + 3)] & (6.12)\cr
\buildrel\rightharpoonup\over{n}_3\cdot
\buildrel\rightharpoonup\over{n}_B\big|_E & = 
2(1-\alpha)^2(1 - 2\alpha)[2(\gamma-1)\alpha -  1] & (6.13) \cr
\buildrel\rightharpoonup\over{n}_3 \cdot \buildrel\rightharpoonup\over{n}_H
\big|_E & = 2(1-\alpha)(1-2\alpha)(1+4\alpha) & (6.14)\cr
}$$
where $\buildrel\rightharpoonup\over{n}_H$ is used to denote an outward normal
to the surface $J^2 = I(1-I)$.  

From now on we restrict ourselves to the case $1< \gamma < 2$, 
since we have done the other cases $\gamma \ge 2$ to earlier papers. 
The integral curves along $ \buildrel\rightharpoonup\over{n}_3$ starting from
$E$ in the direction of increasing $s$  are depicted in Fig. 6.3.

It can be seen  that $\buildrel\rightharpoonup\over{n}_3$ points into different
directions relative to the positions of $A=0$, \ $B=0$, and $H=0$ 
and the sign of the eigenvalue $\lambda_3$ changes when
$\alpha$ varies. Our goal is to use the appropriate direction to 
construct solutions with monotone increasing $s$.
Relative to the surface $B=0$ 
we find that $\buildrel\rightharpoonup\over{n}_3$ points into $\Omega_3$ for
$\alpha\in\left( {1\over 2(\gamma-1)},\ 1\right)$, 
into $\Omega_{33}$ and $B<0$ for $\alpha\in\left({1\over 2},
\ {1\over 2(\gamma-1)}\right)$, and $-\buildrel\rightharpoonup\over{n}_3$
points into $\Omega_3$ for $I\in\left( 0,{1\over 2}\right)$.
Relative to the surface $A=0$ 
we find that $\buildrel\rightharpoonup\over{n}_3$ points 
into $\Omega_{33}$ and $A<0$ for 
$\alpha\in\left({\gamma + 3\over 2(3\gamma -1)}, \ 1\right)$,
into $A>0$ for $\alpha\in\left({1\over 2}, \ {\gamma + 3\over 2(3\gamma -1)}
\right)$, and $-\buildrel\rightharpoonup\over{n}_3$ points 
into $\Omega_3$ for
$\alpha\in\left(0,\  {1\over 2}\right)$. 

Now we can depict the integral curves inside $\Omega_{33}$,  see Figure 6.4. 
\medskip
\midinsert\centerline{\psfig{figure=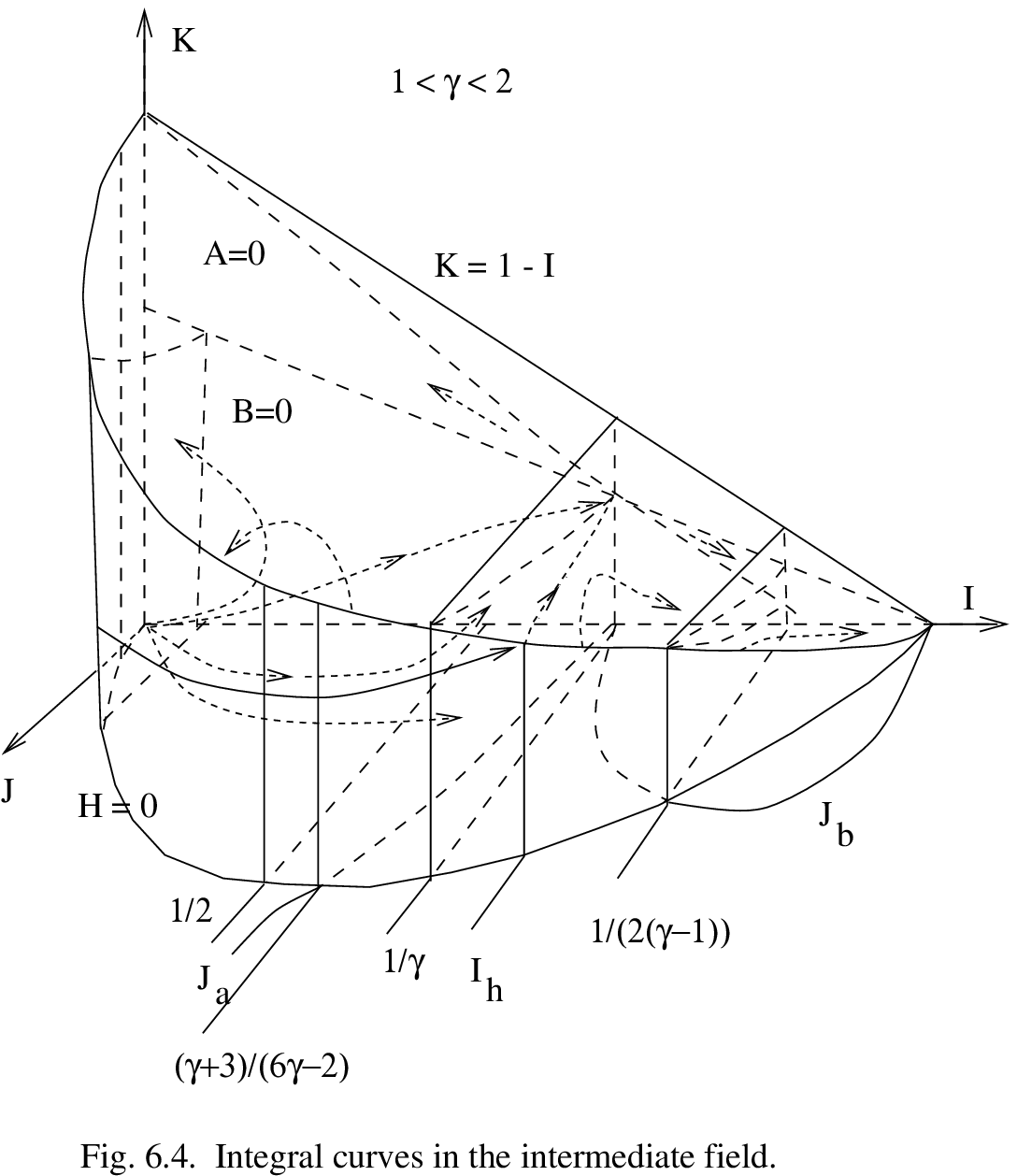,height=5truein}}
\endinsert
\medskip\medskip\medskip\medskip\medskip\medskip\medskip

First we observe that $s$ is an increasing function of $\tau$ inside
$\Omega_{33}$.  $J$ is an increasing function of $\tau$ 
if $0<I<{1\over 2}$, but changes to decreasing when $I\in\left({1\over 2},
\ 1\right)$.  It is useful to observe the stable manifold of the system 
(6.1--3) at
the point (6.6) which contains the transitional integral curve in the case 
$v_0=0$ considered in Sect. 5 and the transitional integral curve from a 
point of $E$ with a critical $I=I_h$. 
There are three kinds of integral curves depending on their positions 
relative to the stable manifold. The first kind 
consists of integral  curves which are below the manifold and 
go to the stationary point $(1,\  0, \  0)$.
Each of the second kind is right on the manifold and goes to the 
stationary point (6.6).  Each of the third kind is above the manifold and
goes to a stationary point on the curve $E$ 
with $I\in \left( 0,\ {1\over 2}\right)$  given in (6.10).  
We mention in particular that no integral curve from inside
$\Omega_{33}$ goes to a point of $E$ between ${1\over 2}<I<1$ because $I$ is an
increasing function of $\tau$ and $\buildrel\rightharpoonup\over{n}_3$ is
pointing towards $(1,0,0)$ for ${1\over 2}<I<1$.  
See Lemma A.1 of the {\it Appendix} for a complete proof. Also there is no
integral curve from inside $\Omega_{33}$ which goes to the point $(0,0,1)$
because $J$ is an increasing function of $\tau$ for $I \in \left( 0,{1\over
2}\right)$.

We need to extend the intermediate field solutions to all $s$.
Integral curves ending at point (6.6) are
already defined for all $s>0$ since $\tau\to +\infty$ and the right-hand
side of equation (6.4) does not vanish at (6.6).
For the first kind of integral curves which end at $(1,0,0)$, we show that $s$
approaches finite values although $\tau\to\infty$.  Since the proof is long and
tedious, we put it in the {\it Appendix} to avoid interuption of our
construction.  We use the natural continuation of vacuum $\rho = 0$ to extend
our solutions till $s = \infty$.  The values of $u,v$ are not needed to be
specified in the vacuum. Near the edge of the vacuum we find 
$$
v = b(r - r^{**})^{\beta}, \qquad \beta \equiv {\gamma+1\over 2(\gamma-1)}
$$
where $b>0$ is any constant. The asymptotics of $u$ and $c$ given in
Sect. 5 are still valid. From these asymptotics we find that $M= \infty$ 
and $M_s = {2\over \sqrt{3-\gamma}}$, which are the same as in Sect. 5.

For the third kind of integral curves which end on the upper half of $E$ with
$\alpha\in \left( 0,{1\over 2}\right)$, 
we show that $s$ approaches finite values also.  
The proof is easy because we have
$$
{ds\over dJ} = {s\over J} {1-I\over 1-2I}\eqno(6.15)
$$
from the second equation of system (4.3).  The right-hand side of (6.15) is
nonsingular for $(I,J,K)\in E$ with $0<I<{1\over 2}$.  Thus $s$ is finite
around any point of $E$ with $0<I<{1\over 2}$.

We comment that integral curves on the surface $J^2 = I(1-I)$ starting from the
origin and ending on a point of $E$ with $\alpha \in (I_h,1)$, where 
$I_h\in\left(1/\gamma, 1/(2(\gamma-1)\right)$,  can be continued
through the direction (6.11) into $\Omega_{33}$, 
and they go toward $(1,0,0)$ with finite ending $s$ values.  
The further extension by vacuum is also valid.  
But integral curves on the surface $J^2 = I(1-I)$ starting from the
origin and ending on a point of $E$ with $\alpha\in(1/2, I_h)$ are continued
along the direction (6.11) into
$\Omega_{33}$ and curl back to a point of $E$ with $0<\alpha< 1/2$.
Finally integral curves on the surface $J^2 = I(1-I)$ starting from the
origin and ending on a point of $E$ with $0<\alpha< 1/2$ will be
continued along direction (6.11) into $H>0$ and discussed in the next section. 
It would be interesting to find explicitly the function $I_h = I_h(\gamma)$. 
We know that $I_h = 1/2$ for $\gamma =2$, see [16, 17].  Numerically, we find
that 
$$
I_h(1.1) = 0.93,\quad I_h(1.4) = 0.77,\quad I_h(1.7)=0.63,
$$
see paper [16] for numerical procedures used. For the special data 
$(u_0,\ v_0, \ \rho_0)$ with $u_0 =0$, the critical value $I_h$ can be 
expressed through the solution formula (6.9) by the critical initial
Mach number
$$
M_h(\gamma) \equiv \sqrt{I_h}/(1-I_h). \eqno(6.16)
$$
Since $I_h \in \left(1/\gamma, 1/(2(\gamma-1))\right)$, we find that
$$
{\sqrt{\gamma}\over \gamma -1} < M_h(\gamma) 
< {\sqrt{2(\gamma-1)}\over 2\gamma-3}. \eqno(6.17)
$$
The lower bound is valid for all $\gamma\in (1, 2)$, but the upper bound
is valid only for $\gamma\in (3/2, 2)$.
So we have constructed all solutions globally except those which end on $E$
with $0<I<{1\over 2}$.
		\bigskip

\noindent {\bf 7.  Inner-field solutions.}
		\medskip
We extend solutions that end on the set $E$  of (6.10) 
with $0<I<{1\over 2}$ in this section.  We
find that these solutions go along the directions
$\buildrel\rightharpoonup\over{n}_3$ given in (6.11) into the region between
$A=0$ and $B=0$ in $0<I<{1\over 2}$ and $J^2>I(1-I)$, and eventually go to
infinity.  The scaled variables $I,J$, and $K$ are not suitable for this 
portion of the solutions.

We restart from system (2.7) with data $(u,v,\rho,r)$ satisfying the 
relations $E:K =
1-I$ and $J^2 = I(1-I)$ for $I\in \left( 0,{1\over 2}\right)$.  
In terms of $(u,v, c, r)$, these data are in the form
$$\cases{
u\big|_{r=\alpha} = \beta \cr
		\noalign{\vskip6truept}
v\big|_{r=\alpha} = \sqrt{\beta(\alpha-\beta)}\cr
		\noalign{\vskip6truept}
c\big|_{r=\alpha} = \alpha-\beta\cr
}\eqno(7.1)
$$
where $\alpha>0$ and $\beta\in \left( 0,{\alpha\over 2}\right)$ are arbitrary.
The directions $\buildrel\rightharpoonup\over{n}_3$ given in (6.11) point into
the region $\Delta >0,\ \Theta >0, \ \Sigma >0$.

Asymptotic analysis shows that problem (2.7) and  (7.1) 
has solutions $(u(r),v(r),\rho(r))$ which vanish as $r\to 0+$. 
The solutions are in the region   $\Delta >0,\ \Theta >0, \ \Sigma >0$.
We perform an asymptotic analysis a
priori to determine the orders at which $(u,v,\rho)$ vanish as $r\to 0+$. 
For polytropic gases $p(\rho) = A_2\rho^\gamma$,
we find 
$$
u={ar\over {\ln{\sigma\over r}}},\quad v={k\over {(\ln{\sigma\over r})^a}},
\quad \rho={d\over {(\ln{\sigma\over r})^{2a}}}\eqno(7.2)
$$
where
$$
a={1\over 2(2-\gamma)}, \qquad k^2={\gamma\over 2-\gamma}d^{\gamma-1},
$$
and $\sigma >0$ and $d>0$ are arbitrary. 
We can easily verify that these directions
are such that $\Delta>0, \Theta>0, \Sigma > 0$ for $r>0$.
In terms of the sound speed $c$, we find
$$
c= {k\sqrt{A_2(2-\gamma)}\over (\ln{\sigma\over r})^{(\gamma-1)a}}.
$$
We are motivated by this asymptotic analysis 
to use the scaled variables
$$
U = {u\over c},\ V = {v\over c}, \ R = {r\over c}\,. \eqno(7.3)
$$
We can rewrite system (2.7) into a new form
$$\cases{
{du\over d\tau} = {r-u\over c^3}\Sigma \cr
		\noalign{\vskip6truept}
{dv\over d\tau} = {uv\over c^3} \Delta \cr
		\noalign{\vskip6truept}
{dr\over d\tau} = {r(r-u)\over c^3}\Delta \cr
		\noalign{\vskip6truept}
{dc\over d\tau} = {\gamma-1\over 2} {r-u\over c^2} \Theta\cr
}\eqno(7.4)
$$
where $\tau$ is a parameter.  In terms of the
variables in (7.3), we find
$$\cases{
{dU\over d\tau} = (R-U)\tilde{\Sigma} - \lambda U(R-U)\tilde{\Theta}
\equiv A_1 \cr
		\noalign{\vskip6truept}
{dV\over d\tau} = UV\tilde{\Delta} - \lambda V(R-U)\tilde{\Theta}
\equiv C_1 \cr
		\noalign{\vskip6truept}
{dR\over d\tau} = R(R-U)\tilde{\Delta} - \lambda R(R-U)\tilde{\Theta}
\equiv B_1 \cr
}\eqno(7.5)
$$
where $\lambda \equiv {\gamma-1\over 2}$ and 
 $\tilde{\Delta} \equiv 1-(U-R)^2,\ \tilde{\Theta} 
\equiv V^2 - U(R-U)$ and 
$\tilde{\Sigma} \equiv (R-U)\tilde{\Theta} - U\tilde{\Delta}$.  
System (7.5) is autonomous for $(U,V,R)$.  
We find that the last equation in (7.4) can be written as
$$
{dc\over cd\tau} = \lambda (R-U)\tilde{\Theta}. \eqno(7.6)
$$
So $c$ can be integrated from (7.6) once $(U,V,R)$ are obtained from (7.5). 
The corresponding data of (7.1) for (7.5) and (7.6) are any stationary point
$(U^*,V^*,R^*, c^*)$ satisfying
$$
R^* - U^* = 1,\ V^{*^2} = U^*,\ 0<2U^*<R^*, \eqno(7.7)
$$
$$
c^*>0.\eqno(7.8)
$$
The asymptotic analysis (7.2) shows that
$$
(U,\ V, \ R) \to (0,\ 0,\ 0), \quad r \to 0+.\eqno(7.9)
$$
After (7.5-9) are solved, we use the third equation in (7.4) to show that $r$
is an increasing function of $\tau\in\Bbb{R}$ and $r\to 0$ or
$\alpha\in(0,+\infty)$ as $\tau$ goes to $\mp \infty$ respectively.

For more refined change of variables, we further find that our asymptotic 
analysis (7.2) implies that
$$
\eqalignno{
U & =a b r \left(\ln {\sigma\over r}\right)^{(\gamma -1)a -1}   \cr
V &= {1\over \sqrt{2-\gamma}}\left(\ln{\sigma\over r}\right)^{-{1\over 2}}
  &(7.10) \cr
R & =b r \left(\ln {\sigma\over r}\right)^{(\gamma -1)a }, \cr
}
$$
where we set $b\equiv (k\sqrt{(2-\gamma)A_2})^{-1}$.
It follows that the above asymptotic behavior yields
the function relations of $(R, U)$ to $V$ 
$$\eqalignno{
R & = {\sigma\over \sqrt{\gamma A_2d^{\gamma -1}}}
\exp\left(-{b^2k^2\over V^2}\right)
\left({bk\over V}\right)^{2(\gamma -1)a} \cr
U & = {1\over 2}RV^2 \cr
}
$$
which are not analytic at $V=0$.
It can be shown that the directions (7.10) enter the region $A_1 >0$, \ 
$B_1>0$, and $C_1>0$. We omit the proof.

We transform the direction $\buildrel\rightharpoonup\over{n}_3$ of (6.11) to
the new variables $(U, \ V, \ R)$ at a point of (7.7). 
It is possible to use
linearization for the derivation, but we prefer a direct transformation. 
Let 
$\buildrel\rightharpoonup\over{n}=( n_1,\ n_2,\ n_3)$ be an arbitrary
vector in the coordinate space $(I,\ J, \ K)$,  which is a tangent 
direction of an  integral curve of (6.1--3). So we have 
$$
{dJ\over dI}  = {n_2\over n_1}, \qquad 
{dK\over dI}  = {n_3\over n_1}.
$$
along the integral curve. Because the variables $(I,\ J, \ K)$ are related to
$(U, \ V, \ R)$
by
$$
I=U/R,\quad J = V/R, \quad K = 1/R,
$$
we find that the differentials $(dR,\ dU, \ dV)$ satisfy
$$
{RdV-VdR\over RdU-UdR}  = {n_2\over n_1},\qquad
{-dR\over RdU-UdR}      = {n_3\over n_1},
$$
which can be solved to yield
$$
{dU\over dR}  = {Un_3-n_1\over Rn_3},\qquad
{dV\over dR}  = {Vn_3-n_2\over Rn_3}.
$$
Thus the corresponding tangent vector in the new coordinate space
$(U, \ V, \ R)$ of $\buildrel\rightharpoonup\over{n}=( n_1,\ n_2,\ n_3)$ 
in the coordinate space $(I,\ J, \ K)$ is a multiple of 
$$
(Un_3-n_1, \ Vn_3-n_2, \ Rn_3), 
\eqno(7.11)
$$
where $( n_1,\ n_2,\ n_3)$ as well as $(U, \ V, \ R)$ are evaluated  at 
the same point $(u,\ v,\ c,\ r)$. For our particular problem 
$\buildrel\rightharpoonup\over{n}_3$ of (6.11), we use $\alpha = U/R$
for a point of (7.7) to find the corresponding direction to be a
multiple of
$$\eqalignno{
\buildrel\rightharpoonup\over{n}_3^* = &
(2+3(3-\gamma)U-5(\gamma-1)U^2,\quad \sqrt{U}[(7-3\gamma)U-(\gamma-1)],\cr
 &\qquad (1+U)[\gamma+3-5(\gamma -1)U]).
&(7.12)\cr}
$$

We determine the relative positions of the regions where
$A_1 > 0, \ B_1> 0,\ $ or $C_1>0$  in the  region $R>U, \ U>0, \ V>0$.
The common intersection $A_1=B_1=C_1=0$ is the curve (7.7) of  stationary
points for all $0<U<R$. The intersection $A_1=B_1=0$ is in the plane 
$R=\gamma U$. The intersection $B_1=C_1=0$ is in the plane $R=2 U$.
 The intersection $A_1=C_1=0$ is in the plane $R={2\over 3-\gamma} U$.
We consider  two domains
$$
\eqalignno{
\Omega_{41}& = \{(U, \ V, \ R)\ | \ U>0, \ V>0, \ R>2U, A_1>0, \ B_1>0 \}\cr
\Omega_{42}& = \{(U, \ V, \ R)\ | \ U>0, \ V>0, \ 2U> R> \gamma U, 
A_1>0, \ C_1>0 \}. \cr
}
$$
See Fig. 7.1.

\medskip
\midinsert\centerline{\psfig{figure=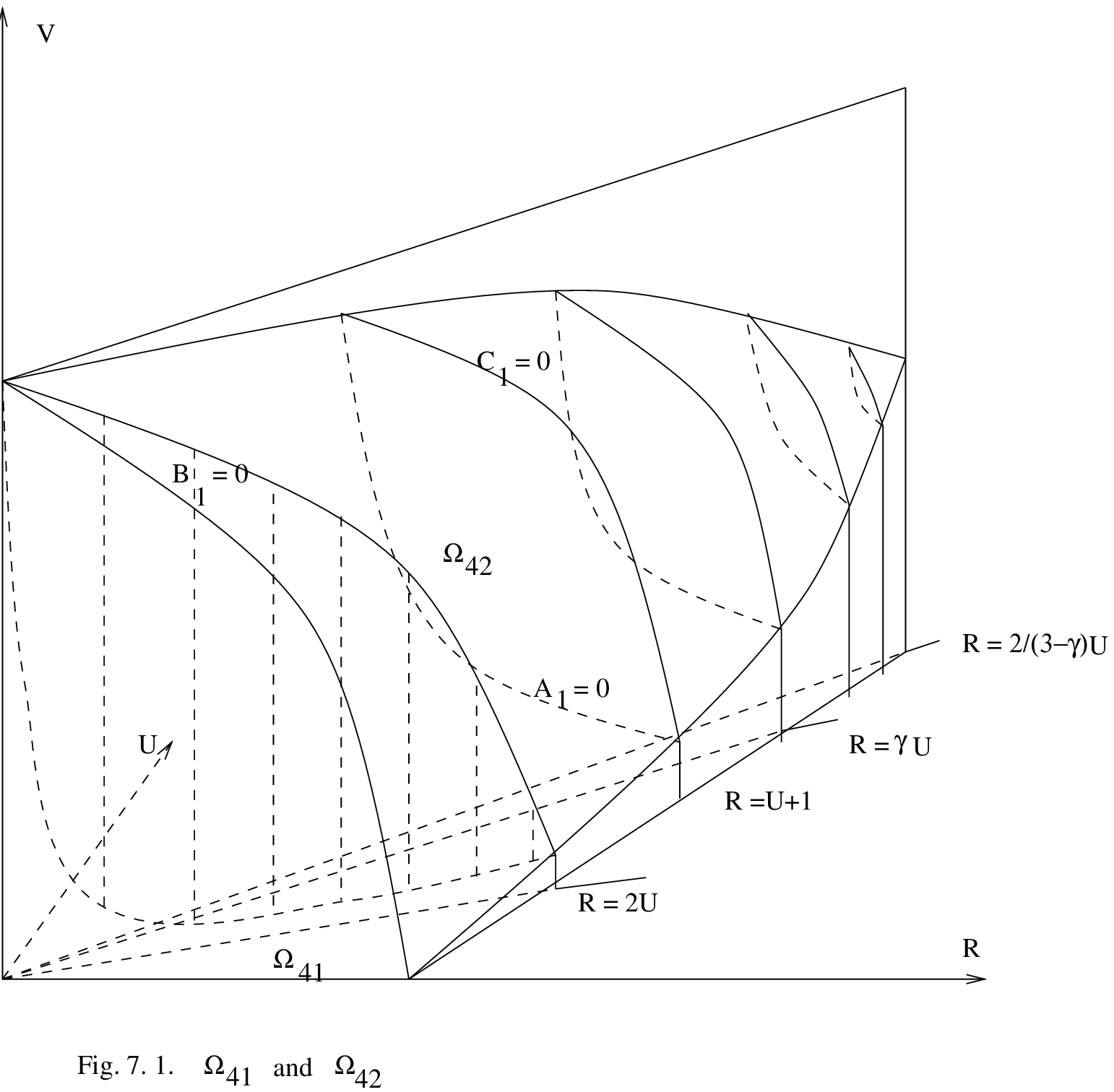,height=6.0truein}}
\endinsert
\medskip\medskip\medskip\medskip\medskip\medskip\medskip

\ 

\

We show that integral curves of (7.5) on the surface $A_1 =0$ of the 
boundaries of both $\Omega_{41}$ and $ \Omega_{42}$ are inward. 
In fact, 
the surface $A_1 =0$ is represented by 
$$
V^2 = {U[2-(\gamma-1)U(R-U)]\over 2R-(\gamma+1)U}
$$
with an outward normal
$$
\buildrel\rightharpoonup\over{n}_{A1} =
\left(n_1, \ -2V, \  {-U[4-(\gamma-1)^2U^2]\over [2R-(\gamma+1)U]^2}\right)
$$
where $n_1$ is some expression that we do not need to compute. We calculate
the inner product of this normal with the tangent direction of integral
curves of (7.5)
$$
{d\over d\tau}(U, \ V, \ R)\cdot\buildrel\rightharpoonup\over{n}_{A1}
 = -2VC_1 -  {U[4-(\gamma-1)^2U^2]\over [2R-(\gamma+1)U]^2}B_1  < 0 
$$
since both $C_1>0$ and $B_1>0$ on $A_1 =0$, in the range $R > \gamma U$.

We show next that integral curves of (7.5) are incoming on the surface
$B_1 =0$ of the domain $\Omega_{41}$ where $R> 2U$. In fact, the surface
$B_1 =0$ is represented by
$$
\lambda V^2 = 1 + \lambda U(R-U)- (R-U)^2
$$
with an inward normal 
$$
\buildrel\rightharpoonup\over{n}_{B1}=
\left( \lambda(R-2U) + 2(R-U), \ -2\lambda V, \ \lambda U -2(R-U)\right).
$$
We find the inner product 
$$
{d\over d\tau}(U, \ V, \ R)\cdot\buildrel\rightharpoonup\over{n}_{B1}
 =A_1[\lambda(R-2U) + 2(R-U)] -2\lambda VC_1 > 0 
$$
since $A_1>0$ and $C_1 < 0$ in the region $R>2U$  on the surface $B_1 =0$.


We show next that integral curves of (7.5) are incoming on the surface
$C_1 =0$ of the domain $\Omega_{42}$ where $2U > R> \gamma U$. 
In fact, the surface $C_1 =0$ is represented by
$$
\lambda V^2 = {U\over R-U} -(1- \lambda)U(R-U)
$$
with an inward normal
$$
\buildrel\rightharpoonup\over{n}_{C1}=
\left({R\over (R-U)^2} +(1-\lambda)(2U-R), \ -2\lambda V,\ 
 -(1-\lambda)U - {U\over (R-U)^2}\right).
$$
We find 
$$
\eqalignno{
{d\over d\tau}(U, \ V, \ R)\cdot\buildrel\rightharpoonup\over{n}_{C1}
& =\left[{R\over (R-U)^2} +(1-\lambda)(2U-R)\right] A_1 -
\left[(1-\lambda)U + {U\over (R-U)^2}\right] B_1 \cr
& > 0 \cr
}
$$
since $A_1>0$, \ $B_1< 0$ on $C_1=0$ in the region ${2\over 3-\gamma}U < R <
2U$ which contains $\gamma U < R < 2U$.

It is easy to see that integral curves of (7.5) are incoming on the surface
$U =0$ of the domain $\Omega_{41}$ since 
$$
{dU\over d\tau} = R^2V^2 > 0
$$
on $U =0$.

We show that integral curves of (7.5) go from $\Omega_{41}$ to $\Omega_{42}$
on the common interface $R=2U$. A normal of $R=2U$ pointing from 
$\Omega_{41}$ to $\Omega_{42}$ is 
$\buildrel\rightharpoonup\over{n} = (2, \ 0, \ -1)$ and
$$
{d\over d\tau}(U, \ V, \ R)\cdot\buildrel\rightharpoonup\over{n}
 = 2U^2(1-U^2-V^2).
$$
on $R=2U$. Using $A_1 \ge 0$ in  $\Omega_{41}$ and  $\Omega_{42}$ on  $R=2U$,
we find
$$
V^2+U^2 - 1  \ge {\lambda +(1-2\lambda)U^2\over 1-\lambda} \ge 0,
$$
thus  integral curves of (7.5) go from $\Omega_{41}$ to $\Omega_{42}$
on the interface $R=2U$.

So integral curves in $\Omega_{41}$ exit $\Omega_{41}$ only through
the side  $R=2U$ or the stationary edge (7.7).

We need to study the local structure of solutions of (7.5) at the 
origin $(U, \ V, \ R)= (0, \ 0, \ 0)$. 
In particular, we would like to prove the asymptotic behavior (7.2).
To reduce the order of degeneracy,
we introduce the variables
$$
X= {U\over R},\quad V_2 = V^2,\quad R_2 = R^2.
$$
We find that
$$
X = a\left( \ln {\sigma\over r}\right)^{-1}, \quad 
V_2 = 2X\eqno(7.13)
$$
as $r\to 0+$ from the asymptotic analysis. 
So the origin  $(X, \ V_2, \ R_2) = (0, \ 0, \ 0)$ will
still be the place for us to look for outgoing solutions. 
Introducing $\tau'$ by 
$$
d\tau' = Rd\tau, \eqno(7.14)
$$
we find the equations (7.5) in terms of  $(X, \ V_2, \ R_2)$ 
to be 
$$\cases{
{dX\over d\tau'}  = (1-X)\left[V_2(1-X)-2X+X(1-X)^2R_2\right]
\equiv (1-X)A_{11} \cr
		\noalign{\vskip12truept}
{dV_2\over d\tau'} =
2V_2 \left\{ X-X(1-X)^2R_2 -\lambda(1-X)\left[V_2 - X(1-X)R_2\right]\right\} 
\cr
		\noalign{\vskip12truept}
{dR_2\over d\tau'} = 
2R_2(1-X)\left\{1-(1-X)^2R_2 - \lambda \left[V_2 - X(1-X)R_2\right]\right\}\cr
}\eqno\eqalign{
(7.15)\cr
(7.16)\cr
(7.17)\cr
}$$
where we have set the expression in the brackets in (7.15) to be $A_{11}$.
We need the local structure of (7.15--17) at $(0, 0, 0)$.
Linearization at the point  $(X, \ V_2, \ R_2)= (0, 0, 0)$ yields
$$\cases{
{dX\over d\tau'}  = - 2X + V_2 \cr
		\noalign{\vskip12truept}
{dV_2\over d\tau'} = 0\cr
	\noalign{\vskip12truept}
{dR_2\over d\tau'} = 2R_2\cr
}
$$
which has eigenvalues
$$
\lambda_1 =-2,\quad \lambda_2 =2,\quad \lambda_3 =0,
$$
and associated eigenvectors
$$
(X, \ V_2, \ R_2)= (1, 0, 0),\  (0, 0, 1), \ (1, 2, 0).
$$

\medskip
\midinsert\centerline{\psfig{figure=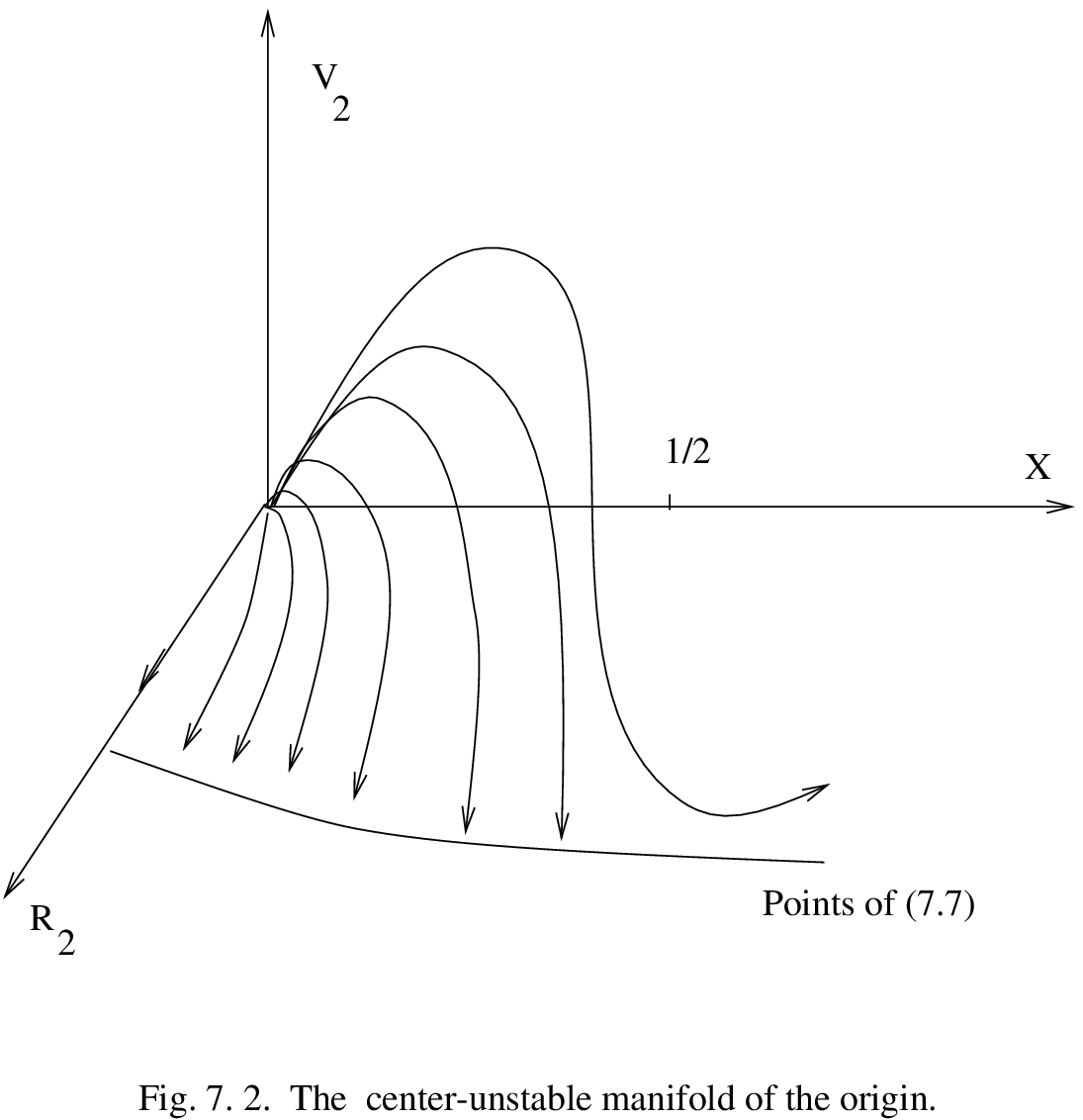,height=4.5truein}}
\endinsert
\medskip\medskip\medskip\medskip\medskip\medskip\medskip

\ 

By the Center Manifold Theorem [8], system (7.15-17) has a $C^k$ 
center manifold which is tangent to the vector $(1, 2, 0)$ for any $k<\infty$. 
Suppose the center manifold takes the form
$$
\eqalignno{
&X  = g(V_2), \qquad g(0)=0, \ g'(0) = 1/2 &(7.18)\cr
&R_2  = h(V_2), \qquad h(0) = h'(0) = 0. &(7.19)\cr
}
$$
Then the flow on the center manifold is described by
$$
{dV_2\over d\tau'} = 2V_2\left\{g-g(1-g)^2h -\lambda(1-g)\left[
V_2-g(1-g)h\right]\right\}.\eqno(7.20)
$$
Using (7.18-19), we can approximate (7.20) to find
$$
{dV_2\over d\tau'} = (2-\gamma)V_2^2 + O(V_2^3).
$$
So the center manifold is unstable since $2-\gamma > 0$.
Now we consider the center-unstable manifold consisting of
the unstable manifold and all center manifolds, see Kelley's paper [8]
for the existence and smoothness. This manifold is two-dimensional. 
It can be shown easily that
they enter the region $A_1 > 0, \ B_1>0, \ C_1>0$.
Furthermore, we can use Henry's approximation of center manifold ([7,8,2])
to find
$$
\eqalign{
g & ={1\over 2}V_2 -{3-\gamma\over 4}V_2^2 + O(V_2^3) \cr
h & =  O(V_2^3), \cr
}
$$
which can be used to show that all center manifolds 
enter the region $A_{11}>0$.

The center-unstable manifold enter the region $\Omega_{41}$. Note that the
unstable manifold, i.e., the $R_2$ axis,  ends on one end of (7.7). We conclude
that nearby integral curves in the center-unstable manifold  end on (7.7) too,
see Fig. 7.2. 

By our earlier analysis, we know that integral curves in 
the center-unstable manifold can only exit $\Omega_{41}$ through the
divider between  $\Omega_{41}$ and  $\Omega_{42}$ or points on (7.7). 
From the continuity of the  center-unstable manifold, we conclude that
every point on (7.7) is covered by the  center-unstable manifold.

Finally, we show that $r\to 0+$ at the origin $(X, V_2, R_2) = (0, 0, 0)$. 
It is easy to see that $c$ is an increasing function of $\tau$ from (7.6)
since $\tilde{\Theta}$ is positive along the center manifolds.  
Thus $c$ is bounded as $\tau\to -\infty$. From (7.3), we have  $r = Rc$.
Hence $r\to 0$ since  $R\to 0$ as $\tau\to -\infty$.

We point out that both the Mach number $M = \sqrt{U^2+V^2}$ and the 
pseudo Mach number $M_s = \sqrt{(U-R)^2 + V^2}$ go to $0$ as $r\to 0$. 
Also, both $U$ and $R-U$ go to $0$ faster than $V$, which explains the
spiralling phenomena in both the physical space and in the self-similar
coordinates.

Thus each point on (7.7) has a solution going to the origin through 
$\Omega_{41}$ with $r\to 0+$.

\bigskip
\noindent {\bf 8.  Conclusions.}
\medskip
We summarize our results.  By a weak solution to the $2-D$ Euler 
equations (2.1)
we mean a bounded vector function $(u,v,\rho)$ satisfying the equations in the
sense of distributions.  Since we deal with only continuous weak solutions in
this paper, we shall not mention any other requirements such as
Rankine-Hugoniot relation or entropy conditions.  let $c = \sqrt{p'(\rho)}$
be the sound speed.  Let $M = \sqrt{u^2+v^2}/c$ be the Mach number. Let
$M_s = [(u-\xi)^2 + (v-\eta)^2]^{1/2}/c$ be the pseudo Mach number of the
self-similar flow in the self-similar coordinates $(\xi, \eta)$. In the
notation of system (2.7), the pseudo Mach number 
$M_s = [(u-r)^2 + v^2]^{1/2}/c$. Set
$$
M_0 \equiv \sqrt{u^2_0 + v^2_0} / {c_0}, \qquad c_0 \equiv \sqrt{p'(\rho_0)}
$$
which are in consistency with previous versions of $M_0$ (Sect. 5 \& 6). Also,
let us introduce the vector parameter (``Mach vector'')
$$
\buildrel\rightharpoonup\over{M}_0 \equiv 
\left( {u_0\over c_0}, {v_0\over c_0}\right).
$$
Recall that
we can assume  $v_0 \ge 0$ without loss of generality. Also, the case $v_0 =0$
has been resolved in [17].
		
\medskip
\midinsert\centerline{\psfig{figure=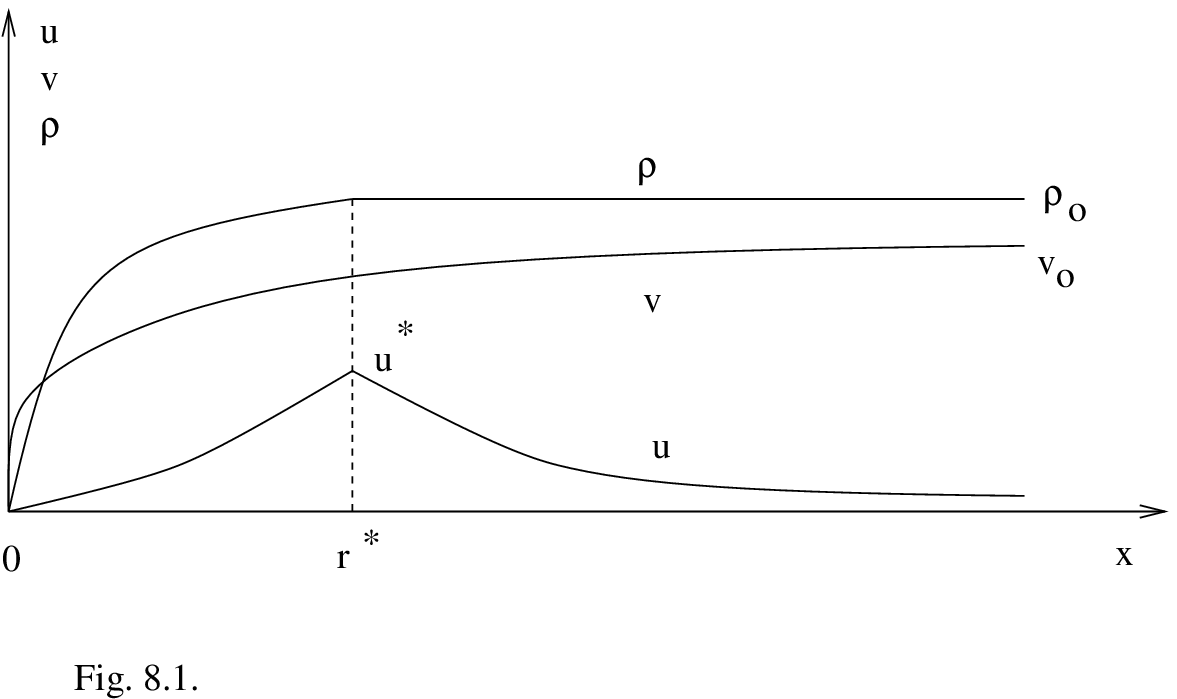,height=2.75truein}}
\endinsert\bigskip\medskip\medskip\medskip\medskip\medskip
\centerline{{\rm \vtop{ 
\hbox{A solution $(u, v, \rho)$ vs. the $x$-axis at time } 
	\vskip-3pt
\hbox{$t=1$ with a datum $u_0=0$ and $M_0 < \sqrt{2}$. }}\hskip1truein }}
\bigskip

\proclaim{Theorem}.  Assume $p(\rho) = A_2\rho^\gamma$ for any constants
 $A_2>0$ and $1<\gamma<2$. Then for any datum $(u_0,v_0,\rho_0)$ 
with $u_0\geq 0$, $v_0 > 0$,  and
$\rho_0>0$, there exists a weak solution $(u,v,\rho)$ to the initial value
problem (2.1) and (2.5). 
The solution is continuous, self-similar, and axisymmetric for $t>0$.  It
takes on initial datum almost everywhere and in $L^q_{{\rm loc}}$ for any $q>1$
as $t\to 0+$.  Additionally, we have the following properties. 
		\medskip
{\rm (I).}  If $u_0=0$, the solution consists of two or three smooth pieces, 
depending on the initial Mach number $M_0$. 
The first piece of the solution is always given explicitly by the formula (3.4)
with constant Mach number near $r=\infty$. 
		\smallskip
{\rm (I.a).} If $M_0 < \sqrt{2}$, the solution consists of two smooth pieces.
The end point of the first piece of the solution is on
the set of stationary points (6.5) with an $I<1/2$, which is transformed to
be the stationary set (7.7) from which the second piece
begins in the direction given in (6.11) or equivalently (7.12). 
It is  an integral curve on 
the center-unstable manifold of the ODE (7.15-17) at the origin.
See Fig. 7.2 and 8.1. Both the Mach and the pseudo Mach numbers are zero
at the center point of the solution.

\smallskip
{\rm (I.b).} If $M_0 > \sqrt{2}$, but $M_0 < M_h(\gamma)$ 
(which is defined in (6.16) ), the solution consists of three 
smooth pieces. The second piece is given by the solution of the 
ODE(4.3) with direction (6.11) such that the solution goes to a stationary
point on (6.5) with an $I<1/2$. The third piece is identical to 
the second piece of a solution described in the previous case $M_0 < \sqrt{2}$.
\smallskip
{\rm (I.c).} If $M_0 = M_h(\gamma)$, the solution consists of two pieces and 
the second piece is given by the solution of the 
ODE(4.3) with direction (6.11) such that the solution goes to the stationary
point (6.6).
		\smallskip
{\rm (I.d).} If $M_0>M_h(\gamma)$, the solution consists 
of three  smooth pieces and the second piece is given by the ODE (4.3) with 
direction  (6.11), see Fig. 8.2.   The radial and pure
rotational conponents and the density function of the solution are increasing
functions of the spatial radius $r$ in the region $u^{**}t<r<r^*t$ where $r^*$
is given by (3.5) and $u^{**}>0$ is the radial velocity at the inner 
end of the second piece.  The third piece is the vacuum $\rho=0$ with domain 
$0\leq r < u^{**}t$. At the edge of the vacuum, there is no rotation,
the Mach number is infinity, and the pseudo Mach number 
$M_s ={2\over \sqrt{3-\gamma}}$ (for $\gamma <3$).

\smallskip
{\rm (I.e).} The critical Mach number $M_h(\gamma)$ has a lower bound
$$
\sqrt{\gamma}/(\gamma -1) < M_h(\gamma)$$
for all $\gamma \in (1, 2)$. It has an upper bound
$$
M_h(\gamma) < \sqrt{2(\gamma-1)}/(2\gamma -3)
$$
if $\gamma \in (3/2, 2)$.
		\medskip
{\rm (II).}  Assume  $u_0>0,\ \rho_0>0$, and $v_0> 0$.  
Depending on the vector parameter $\buildrel\rightharpoonup\over{M}_0$, 
the solution can be globally smooth, contain a region of vacuum as in 
case (I.a), or contain an inner piece of a solution of (I.d). 
In terms of the selfsimilar coordinates and variables, the ODE (6.1-3)
has at the point (6.6) a stable manifold, which contains the transitional
solutions in the special cases $u_0 = 0$ of case (I) and $v_0 = 0$. If 
the vector parameter $\buildrel\rightharpoonup\over{M}_0$
is such that the solution is on the stable
manifold, then the entire solution consists of one smooth piece and is 
governed by the system of ODE (3.1). 
If the vector parameter $\buildrel\rightharpoonup\over{M}_0$ 
is such that the solution is below the stable
manifold, then the solution consists of a vacuum region at the center and 
the far-field solution given by the system of ODE (3.1) with data (3.2).  
If the vector parameter $\buildrel\rightharpoonup\over{M}_0$ 
is such that the solution is above the 
stable manifold, then the solution consists of two pieces, one is given
by the system of ODE (3.1) with data (3.2), the other by an integral curve
in the center-unstable manifold near the origin of (7.15-17), 
and the density $\rho$ of the solution vanishes only at one point, the origin
of space.
\medskip

We omit the proof of the theorem with the remark that 
it can be verified that the pieces of integral curves when
combined continuously in the order specified in the theorem form weak 
solutions to (2.1)(2.4). 

{\bf The case $\gamma =1$}. We now point out necessary changes in Sections
4--7 to accommodate the case $\gamma =1$. The equations (4.3) are still 
valid although the equation for $K$ is trivial since the sound speed $c$
is a constant. We need to supplement (4.3) with the equation for $\rho$
in (3.1), the equation for $\rho$ being decoupled from (4.3) now. 
Conclusions in Sect. 5 are modified as follows: all integral curves in 
$\Omega$ go to the point $(I, K) = (0, 1)$ except the trivial one 
$K=0$. In Sect. 6 all integral curves in $\Omega_3 =\Omega_{33}$
from the origin $(I, J, K)= (0, 0, 0)$ enter the set $E$ with $\alpha < 1/2$.
In Sect. 7 the conclusion is the same, but the proof, which we shall omit,
needs nontrivial changes since the domains $\Omega_{41}$ and $\Omega_{42}$
are not bounded anymore. The overall conclusion for $\gamma =1$ is the 
same as stated in the previous theorem, the main difference  is that
the threshold Mach number $M_h(1) = \infty$ and the threshold stable
manifold mentioned in Step II is the bottom surface $K=0$, thereby all 
solutions in the case $\gamma =1$ have point cavities at the centers.


\medskip
\midinsert\centerline{\psfig{figure=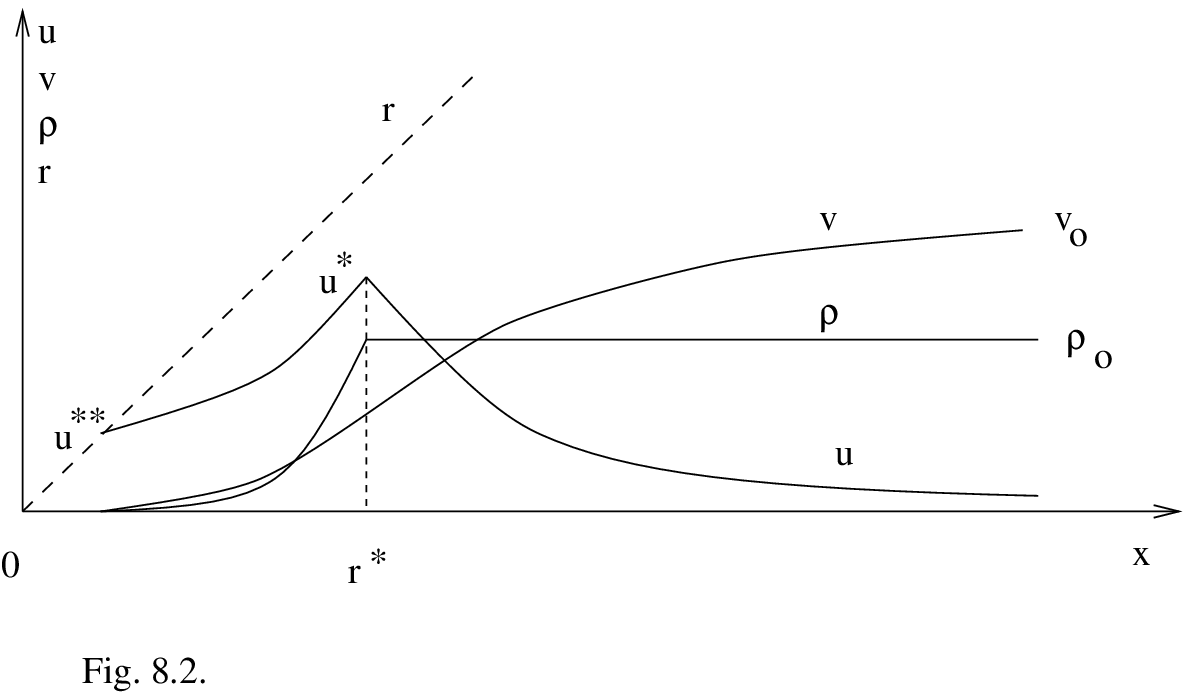,height=2.75truein}}
\endinsert\medskip\medskip\medskip\medskip\medskip\medskip\medskip
\centerline{{\rm \vtop{ 
\hbox{A solution $(u, v, \rho)$ vs. the $x$-axis at time } 
	\vskip-3pt
\hbox{$t=1$ with a datum $u_0=0$ and $M_0 > M_h(\gamma)$. }}\hskip1truein }}

\vfill\eject

\noindent{\bf Appendix.}  Finiteness of the parameter $s$ at point $(1,0,0)$

In this appendix we show mainly that the parameter $s$,
the reciprocal of the radius $r$, approaches a finite 
value as an integral curve approaches the stationary point $(1,0,0)$ in 
$\Omega_3$ defined in section 6.
		\bigskip

\proclaim{Lemma A.1}.  Assume $\gamma>1$. Then for any integral curve of 
the  system (6.1-4)
that goes from $I < 1/2$ to $I>1/2$  inside the domain $\Omega_{33}$, 
there exists an $\varepsilon\in (0,1)$ such that the integral curve is 
inside the cylinder
$$
C_\varepsilon\equiv J^2 - \varepsilon I(1-I)<0 \eqno(A.1)
$$
whenever its $I$ component is in $\left({1\over 2},1\right)$.
		\bigskip

\noindent {\bf Proof.}  We choose $\varepsilon\in (0,1)$ such that the integral
curve is inside the cylinder $J^2 < \varepsilon I(1-I)$ at the special point $I
= {1\over 2}$.  For this $\varepsilon$, we show that integral curves on the
surface of the cylinder $J^2 = \varepsilon I(1-I)$ inside $\Omega_{33}$ in the
portion
$I\in
\left( {1\over 2},1\right)$ are all going into the cylinder $J^2<\varepsilon
I(1-I)$.  We first calculate an outward normal to the cylinder $J^2 =
\varepsilon I(1-I)$ to be
$$
\buildrel\rightharpoonup\over{n}_\varepsilon = (-\varepsilon(1-2I), 2J, 0).
$$
We then calculate the inner product of the normal
$\buildrel\rightharpoonup\over{n}_\varepsilon$ with the tangent direction of
any integral curve on the cylinder $J^2 = \varepsilon I(1-I)$
$$\eqalign{
\buildrel\rightharpoonup\over{n}_\varepsilon \cdot {d\over d\tau}
(I,J,K)\big|_{J^2 = \varepsilon I(1-I)} & =
-\varepsilon(1-2I)(1-I)[I(1-I)^2(1+\varepsilon) - 2IK^2] \cr
&\qquad + 2\varepsilon I(1-I)(1-2I)\cdot [(1-I)^2-K^2] \cr
& = \varepsilon I(1-2I)(1-I) \{ - (1+\varepsilon)(1-I)^2 + 2K^2 + 2(1-I)^2 -
2K^2\} \cr
& = \varepsilon (1-\varepsilon) I(1-I)^3(1-2I) < 0 \ \hbox{for}\ I\in \left(
{1\over 2},1\right).\cr
}$$
Hence the chosen integral curve will remain inside the cylinder
$J^2<\varepsilon I(1-I)$ in $I\in \left( {1\over 2},1\right)$.  This completes
the proof of Lemma A.1.
		\smallskip
\noindent{\it Remark.} From the proof of Lemma A.1 it is clear that any
integral curve entering $\Omega_{33}$  from the edge $E$ (6.10) will
get into a cylinder $J^2<\varepsilon I(1-I)$ for some $\varepsilon \in (0,1)$
for $I$ close to 1. 

\bigskip
\proclaim{Lemma A.2}. Assume $\gamma \ge 3/2$. Then for any integral curve 
that ends at $(1,0,0)$ of the system (6.1-4) in the domain $\Omega_3$, 
there are three numbers $\varepsilon \in (0,1),\ \beta \in (2,\infty)$, 
and $\tilde{I}\in \left( {1\over 2},1\right)$
such that the integral curve is inside the cylinder $C_\varepsilon <0$ and
below the surface
$$
B_\beta \equiv (1-I)J^2 + I(1-I)^2 - \beta IK^2 = 0 \eqno(A.2)
$$
for all $I \in (\tilde{I},1)$.
		\bigskip

\noindent {\bf Proof.}  We first compute an outward normal
$\buildrel\rightharpoonup\over{n}_\beta$ to the surface $B_\beta = 0$.  
$$
\buildrel\rightharpoonup\over{n}_\beta = 
(\beta K^2 - (1-I)^2 + 2I(1-I) + J^2,\
-2(1-I)J,\ 2\beta IK).
$$
We calculate the inner product of this $\buildrel\rightharpoonup\over{n}_\beta$
with the tangent vector of any integral curve on the surface $B_\beta=0$:
$$\eqalignno{
& \buildrel\rightharpoonup\over{n}_\beta \cdot {d\over
                                        d\tau}(I,J,K)\big|_{B_{\beta=0}} \cr
= & \left[\beta K^2-(1-I)^2 + 2I(1-I)+J^2\right] \cdot
     (1-I) \left[(1-I)^2I + (1-I)^2 - 2IK^2\right] \cr
&\qquad - 2(1-I)J^2(1-2I)[(1-I)^2-K^2] +\cr
&\qquad \beta IK^2(1-I) 
       \left[(\gamma-1)J^2 - (\gamma-1)I(1-I) + 2(1-I)^2 - 2K^2\right] \cr
 = &\left[ 2I(1-I)+{1\over I} J^2 \right] (1-I)^2 \left( 1 - {2\over
\beta}\right) [J^2 + I(1-I)] \cr
& + 2(1-I)^2(1-2I)J^2 {1\over \beta I} \left[J^2 - (\beta-1)I(1-I)\right] 
      + (1-I)^2 \left[J^2 + I(1-I)\right] {1\over I}\cdot \cr
&\Bigg\{\left[(\gamma-1)I -{2\over \beta}(1-I)\right]J^2 -(\gamma-1)I^2(1-I) 
 + 2\left( 1-{1\over \beta}\right) I(1-I)^2 \Bigg\}\,. \cr
}$$
Rewriting the last expression in the form of a polynomial of $J$, we obtain 
$$\eqalignno{
& \buildrel\rightharpoonup\over{n}_\beta \cdot {d\over d\tau}
        (I,J,K)\big|_{B_{\beta=0}} \cr
 = & J^4 \Big\{ {1\over I}(1-I)^2 \left( 1-{2\over \beta}\right) 
     + {2\over \beta I} (1-I)^2(1-2I) + {(1-I)^2\over I} \left[
         (\gamma-1)I - {2\over \beta}(1-I)\right]\Big\}        \cr
& + J^2 \Big\{ 2I(1-I)^3\left( 1-{2\over \beta}\right) + (1-I)^3 
\left( 1-{2\over \beta}\right) - 2 {\beta-1\over \beta} (1-I)^3(1-2I) \cr
& + (1-I)^3 \left[ (\gamma-1)I-{2\over \beta}(1-I)\right] 
+ 2\left( 1-{1\over \beta}\right) (1-I)^4 - (\gamma-1)I(1-I)^3 \Big\} \cr
&  + 2I^2 (1-I)^4 \left( 1-{2\over \beta}\right)- (\gamma-1)I^2(1-I)^4 + 2
\left( 1-{1\over \beta}\right) I(1-I)^5 \cr
 = & J^4 {(1-I)^2\over \beta I} \{\beta-2 + [\beta(\gamma-1)-2]I\} + 
      + J^2 (1-I)^3 \left[ 4I\left( 1-{1\over \beta}\right) + 1- {4\over \beta}
         \right] \cr
& + I(1-I)^4 \left[ 2I\left( 1-{2\over \beta}\right) - (\gamma-1) I + 2
\left( 1-{1\over \beta}\right)(1-I) \right]\,.\cr
 }$$
Observing that both coefficients of $J^4$ and $J^2$ are positive 
when $\beta>2,\ \gamma\geq 2$, and $I \in \left( {1\over 2},1\right)$, we find
in this case in the region
$C_\varepsilon < 0$ that 
$$\eqalignno{
\buildrel\rightharpoonup\over{n}_\beta\cdot {d\over d\tau}
(I,J,K) & \big|_{B_{\beta=0},C_\varepsilon<0} < {\varepsilon^2\over \beta}
I(1-I)^4 \{\beta-2+[\beta(\gamma-1)-2]I\} \cr
&\qquad + \varepsilon I(1-I)^4 \left[ 4I\left( 1-{1\over \beta}\right) + 1 -
{4\over \beta}\right] \cr
&\qquad + I(1-I)^4 \left[ -{2I\over \beta} - (\gamma-1)I + 2\left( 1-{1\over
\beta}\right)\right]  &(A.3) \cr
& = I(1-I)^4 \Big\{ \varepsilon^2 \left[ 1-{2\over \beta} + \left(
\gamma-1-{2\over \beta}\right)I\right] + 4\varepsilon I\left( 1-{1\over
\beta}\right) + \cr
& \qquad  +\varepsilon-{4\varepsilon\over \beta} - \left( {2\over \beta} 
+ \gamma-1\right) I + 2\left( 1-{1\over \beta}\right)\Big\} \cr
& \equiv I(1-I)^4 F(\varepsilon,\beta,I,\gamma) \cr
}$$
where $F$ denotes the expression in the curly brackets. 

We show that the inequality (A.3) still holds even when $\gamma \in [3/2, 2)$.
The quadratic polynomial of $J^2$ attains its maximum at
$$
(J^2)_{max} = {I(1-I)[4I(\beta -1)+\beta -4]\over 
2I[2-\beta(\gamma -1)] +4-2\beta}.
\eqno(A.4)
$$
At $\beta =2$, we calculate
$$
(J^2)_{max} = {1-I\over 2(2-\gamma)} +O((1-I)^2)\eqno(A.5)
$$
for $I$ close to 1.
Since $(J^2)_{max}$ is a continuous function of $\beta$ near 2, we see that
$(J^2)_{max}$ is positive for all $\gamma \in (1,2)$ 
when $\beta$ is close to 2 and $I<1$ is close to 1.
Hence the quadratic polynomial of $J^2$ is an increasing function of $J^2$
in the interval $[0, (J^2)_{max}]$. We show next that
$$
\varepsilon I(1-I) \le (J^2)_{max} \eqno(A.6)
$$
for all $\gamma \in [3/2, 2)$, and $\varepsilon <1$ when $I$ is close to 1
and $\beta$ is close to 2.
In fact, (A.6) holds by continuity, the fact (A.5) and the inequality
$$
\varepsilon < {1\over 2(2-\gamma)}
$$
which holds for all $\gamma \in [3/2, 2)$. So inequality (A.3) still holds 
in the region $C_\varepsilon < 0 $ even when $\gamma \in [3/2, 2)$. 

We need $F$ to be negative near $I = 1$ for some $\varepsilon>0,\ \beta>2$, and
$\gamma\ge 3/2$. For any fixed integral curve that goes to the point $(1,0,0)$,
we first choose
$\varepsilon\in (0,1)$ so that the integral curve lies inside $C_\varepsilon
<0$ for $I \in (\tilde{I}, 1)$ for some $\tilde{I} < 1$.  
Fix this $\varepsilon$.  We calculate the value of
$F$ at the extreme point $\beta=2, I=1$:
$$\eqalignno{
F(\varepsilon,2, 1,\gamma) & =\varepsilon^2(\gamma-2)+\varepsilon -\gamma +1
 \cr
& = (1 - \varepsilon)[1+ 2\varepsilon - \gamma(1+\varepsilon)] \cr 
&\le (1 - \varepsilon)[1+ 2\varepsilon - {3\over 2}(1+\varepsilon)] \cr
&= - {1\over 2}(1 - \varepsilon)^2 < 0. \cr
 }
$$
Since $F$ is a continuous function of $(\beta,I)$ near the point $(2,1)$, we
conclude that for each $\gamma \ge 3/2$ there exists an 
$\tilde{I}\in \left( {1\over 2},1\right)$ such that $F(\varepsilon,\beta,I,
\gamma)<0$ for all $I\in [\tilde{I},1)$ when $\beta$ is
close to $2$.  Since $B_\beta=0$ is close to $B=0$, we can choose a $\beta>2$
such that the integral curve lies below $B_\beta=0$ at $I = \tilde{I}$.  This
integral curve will remain under the surface $B_\beta=0$ for all $I>\tilde{I}$
because of the sign $F<0$.  This completes the proof of Lemma A.2.
		\bigskip

Lemma A.2 does not hold for $\gamma < 3/2$. We need to deal with the case 
$\gamma < 3/2$ separately. Observe that the intersection
point of the  curve $J_b=0$  with the stationary edge (6.5) is topologically
different for $\gamma < 3/2$, see Fig.6.1. 
\smallskip

\proclaim{Lemma A.3}.  Assume $1< \gamma < 3/2$. Then integral curves in 
$\Omega_3$ of system (6.1--3) on the surface of the domain
$$
(2-\gamma)J^2 \le \varepsilon(\gamma I -1)(1-I), \eqno (A.7)
$$
where $\varepsilon\in(0,1)$, always enter the domain.
\smallskip
\noindent {\bf Proof.} We calculate an outward normal to (A.7):
$$
\buildrel\rightharpoonup\over{n} = (- \varepsilon(\gamma + 1 -2\gamma I), \  
2(2-\gamma)J,\ 0).
$$
We then calculate the inner product of this normal with the tangent vector
of an integral curve of (6.1--3) on the surface of the domain (A.7):
$$
\eqalignno{
\buildrel\rightharpoonup\over{n}\cdot{d\over d\tau}(I, J, K) &=
\varepsilon(1-I)(2\gamma I - \gamma - 1)B - 2(2-\gamma)J^2(2I-1)\left[
(1-I)^2 - K^2\right]     &(A.8)\cr
& = \varepsilon(1-I)^2(2\gamma I - \gamma - 1)\left[J^2+I(1-I)\right] 
- 2(2-\gamma)J^2(2I-1)(1-I)^2  \cr
&\qquad\qquad + K^2\left[2(2-\gamma)J^2(2I-1)-
2\varepsilon I(1-I)(2\gamma I - \gamma - 1)\right]\cr
& = \varepsilon(1-I)^2(2\gamma I - \gamma - 1)[J^2+I(1-I)] 
- 2(2-\gamma)J^2(2I-1)(1-I)^2 \cr
&\qquad\qquad + 2\varepsilon(1-I)^2K^2. \cr
}
$$
We observe that this inner product is an increasing function of $K^2$.
In the domain $\Omega_3$, the largest $K^2$ is achieved on the surface
$B=0$. From the step (A.8), we find that this inner product in negative
for $I>1/2$. This completes the proof of Lemma A.3.
\bigskip

\proclaim{Lemma A.4}. Assume $1<\gamma<3/2$. Then for any integral curve 
that gets into (at least once) the domain 
$$
(2-\gamma)J^2 \le (\gamma I -1)(1-I) \eqno(A.9)
$$
and ends at $(1,0,0)$ of the system (6.1-4) in the domain $\Omega_3$, 
there are three numbers $\varepsilon \in (0,1),\ \beta \in (2,\infty)$, 
and $\tilde{I}\in \left( {1\over 2},1\right)$
such that the integral curve is inside the cylinder (A.7) and
below the surface (A.2).
$$
B_\beta \equiv (1-I)J^2 + I(1-I)^2 - \beta IK^2 = 0 \eqno(A.2)
$$
when $I \in (\tilde{I},1)$.
\smallskip

\noindent {\bf Proof.} The proof parallels that of Lemma 2: use the surface
of the domain (A.7) in places where $C_\varepsilon = 0 $ is used. 

\bigskip

\proclaim{Theorem A}.  For all $\gamma>1$, the parameter $s$ is finite 
for any integral curve of
(6.1-4) that goes to the point $(1,0,0)$ inside $\Omega_3$.
		\bigskip

\noindent {\bf Proof.}  (i) Assume $\gamma \ge 3/2$ or the integral curve
gets into the domain (A.7) at least once, hence the integral curve 
remains in $B_\beta>0$ for $I$ close to $1$.
We find from the first equation of (4.3) that
$$
s{dI\over ds} = 2I - (1-I) {I(1-I)-J^2\over (1-I)^2-K^2} \eqno(A.10)
$$
The right hand side of (A.10) is a decreasing function of $K^2$ in $\Omega_3$. 
So for $I$ close to
$1$ we find
$$\eqalignno{
s{dI\over ds} & \geq 2I-(1-I) I {I(1-I)-J^2 \over 
\left( 1-{1\over \beta}\right)
I(1-I)^2 - {1\over \beta}(1-I)J^2} \cr
& = 2I - I\beta {I(1-I)-J^2\over (\beta-1)I(1-I)-J^2} = 2I - \beta I +
{\beta(\beta-2)I^2(1-I)\over (\beta-1)I(1-I)-J^2}\,. & (A.11)\cr
}$$
The very last expression of (A.11) is an increasing function of $J^2$, so we
further find by using $J^2\geq 0$ that
$$
s {dI\over ds} \geq (2-\beta)I + {\beta(\beta-2)\over \beta-1}I = {\beta-2\over
\beta-1}I.
$$
So it can only take a finite amount of $s$ for $I$ to reach $1$.
(ii) If $1<\gamma<3/2$ and the integral curve in entirely outside the 
domain (A.9) when 
$I$ is close to 1, we have 
$$
(2-\gamma)J^2 \ge (\gamma I -1)(1-I) \eqno(A.12)
$$
for $I$ close to 1. We can use equations (6.2) and (6.4) to find
$$
{dJ\over ds} = {J(1-2I)\over s(1-I)}. \eqno(A.13)
$$
Using (A.12) in (A.13), we find
$$
{d\ln ({1\over J})\over d\ln s} = {2I-1\over 1-I} > \alpha^2({1\over J})^2
$$
where $ \alpha^2$ is a positive constant. It is easy to derive that 
$s<\infty$ from this last differential inequality.

	\vskip.5truein
\noindent E-mail addresses: yzheng@indiana.edu, tzhang@math03.math.ac.cn

\noindent Acknowledgement --

We appreciate helpful discussions with Taiping Liu and Zhihong (Jeff) Xia.  
This work is supported by NSF DMS--9303414 and  Alfred P. Sloan Research 
Fellows award for Yuxi Zheng and by National Fundamental
Research Program of State Commission of Science and Technology of China 
and NSFC for Tong Zhang.

\vfill\eject

\noindent{\bf Figure Captions.}
\bigskip

Fig. 1.1. The velocity vector of a typical swirling solution. 
Here $\gamma =1.4$, $u_0 =0.0$, $v_0 =1.0$, $\rho_0=1.0$, and $A_2=1.0$. 
The figure shows the solution at time $t=1$ in the square $|x|\le 1.7535$,
and $|y|\le 1.7535$. 

Fig. 1.2. The fluid density of the typical swirling solution in Fig. 1.1.
Outside the square the density is the constant $1$.

Fig. 5.1. Phase portrait of solutions without swirls.

Fig. 5.2. Estimate of the transitional Mach number $M_h(\gamma) <
{\sqrt{2}\over \gamma-1}$.  In this figure $d =(\gamma -1)^2/2$, 
$f=M_h^{-2}(\gamma)$, and $N$ and $D$ are where the numerator
and denominator of the right-hand side of (5.8) vanish respectively.

Fig. 5.3. A zero-swirl transitional solution.

Fig. 5.4. The parameter $s$ is finite near $(1,0)$.

Fig. 6.1.  The region $\Omega_3$ and $\Omega_{33}$ .

Fig. 6.2.  Local structure of solutions at the stationary point (6.6).

Fig. 6.3.  Solutions in the intermediate field for data with $u_0=0$.

Fig. 6.4.  Integral curves in the intermediate field.

Fig. 7.1.  The domains $\Omega_{41}$ and  $\Omega_{42}$.

Fig. 7.2.  The center-unstable manifold. 

Fig. 8.1.  A solution $(u, v, \rho)$ vs. the $x$-axis at time $t=1$ with
            a datum $u_0=0$ and $M_0 < \sqrt{2}$.

Fig. 8.2.  A solution $(u, v, \rho)$ vs. the $x$-axis at time $t=1$ with
            a datum $u_0=0$ and $M_0 > M_h(\gamma)$.

\vfill\eject
\centerline{Reference}
		\medskip
\item{[1]} { Bellamy-Knights, P. G.}: Viscous compressible heat conducting
spiralling flow, Q. J. Mech. Appl. Math. {\bf 33} (1980), Pt. 3.
 
\item{[2]} { Carr, J.}: Applications of centre manifold theory, 
Applied Mathematical Sciences, 35, Springer-Verlag, New York, 
Heidelberg, Berlin. 1981.

\item{[3]} { Colonius, T., Lele, S. K., Moin, P.}: The free compressible
viscous vortex, J. Fluid Mech., {\bf 230} (1991), 45--73.
 
\item{[4]} { Courant, R., Friedrichs, K. O.}: Supersonic flow and 
shock waves,
Applied Mathematical Sciences 21, Springer-Verlag, 1948, 1976.
 
\item{[5]} { Chen, G. -Q., Glimm, J.}:  Global solutions to the compressible
Euler equations with geometrical structure, preprint (1994).
 
\item{[6]} {Guckenheimer, J., Holmes, P.}: Nonlinear oscillations,
dynamical systems, and bifurcations of vector fields, Applied 
 Mathematical Sciences, 42, Springer-Verlag, New York, 
Heidelberg, Berlin. 1983.

\item{[7]} { Henry, D.}: Geometric theory of semilinear parabolic Equations,
Springer Lectures Notes in Mathematics, Vol. 840. Springer-Verlag, New York, 
Heidelberg, Berlin. 1981.

\item{[8]} { Kelley, A.}: The stable, center stable, center, 
center unstable and unstable manifolds, J. Diff. Eqns., 
{\bf 3} (1967), 546--570

\item{[9]} { Mack, L. M.}: The compressible viscous heat-conducting vortex,
J. Fluid Mech., {\bf 8} (1960), 284--292.
 
\item{[10]} { Majda, A.}: Vorticity, turbulence, and acoustics in fluid flow,
SIAM Review, {\bf 33} (1991), 349--388.

\item{[11]} { Powell}: Lecture series at the von Karman Institute, 1990
 
\item{[12]} { Serrin, J.}: The swirling vortex, Phil. Trans. Roy. Soc. 
London, A., {\bf 271} (1972), 325--360.

\item{[13]}  { Dechun Tan and Tong Zhang}:  Two-Dimensional Riemann Problem 
for a 2x2 Hyperbolic System of Nonlinear Conservation Laws (I) and (II),
J. Diff. Equations,  {\bf 111}(1994), 203--282.

\item{[14]} { Wagner , D.}:  The Riemann Problem in Two Space Dimensions 
for a Single Conservation Law, SIAM  Journal on  Mathematical Analysis, 
{\bf 14} (1983), 534--559

\item{[15]} { Zhang, Tong and Zheng, Yuxi}:  Conjecture on Structure of 
Solutions of Riemann Problem for 2-D Gasdynamic Systems, 
SIAM  Journal on  Mathematical  Analysis, {\bf 21} (1990), 593--630.

\item{[16]} {  -----------}: Exact Spiral Solutions of the Two Dimensional 
Compressible Euler Equations,  preprint of the Institute for Scientific
Computing and  Applied Mathematics, No. 9410, 1994. Submitted to
Discrete and Continuous Dynamical Systems.

\item{[17]} { Zheng, Yuxi and Zhang, Tong}:  
Axisymmetric Solutions of the Euler Equations 
 for Square Polytropic Gases, preprint, Dept. of Math., Indiana University.
Submitted to Arch. Rat. Mech. Anal. in April 1996.

\end